\documentclass[11pt]{amsart}
\usepackage{amscd, amsfonts, amsthm, graphicx, amssymb}
\title{Groups of volume-preserving diffeomorphisms of noncompact manifolds and mass flow toward ends} 
\author{Tatsuhiko Yagasaki}
\address{Division of Mathematics,  Graduate School of Science and Technology, 
Kyoto Institute of Technology, Kyoto, 606-8585, Japan}
\email{yagasaki@kit.ac.jp}

\newtheorem{theorem}{Theorem}[section]
\newtheorem*{theorem-v-1}{Theorem 1.1$'$}
\newtheorem*{theorem-a-1}{Theorem 1.2$'$}
\newtheorem{prop}{Proposition}[section] 
\newtheorem{cor}{Corollary}[section] 
\newtheorem{lemma}{Lemma}[section]

\theoremstyle{definition}
\newtheorem{defi}{Definition}[section]
\newtheorem{remark}{Remark}[section]
\newtheorem{notation}{Notation}[section]

\newtheorem{assumption}{Assumption}[section]

\def \cal {\mathcal}
\def \phi {\varphi}
\def \ds {\displaystyle}

\def \lra {\longrightarrow}

\def \e {\varepsilon}

\begin{document}
\baselineskip 6 mm

\thispagestyle{empty}

\maketitle

\begin{abstract} 
Suppose $M$ is a noncompact connected oriented $C^\infty$ $n$-manifold and 
$\omega$ is a positive volume form on $M$. 
Let ${\cal D}^+(M)$ denote the group of orientation preserving diffeomorphisms of $M$ 
endowed with the compact-open $C^\infty$ topology and 
${\cal D}(M; \omega)$ denote the subgroup of $\omega$-preserving diffeomorphisms of $M$. 
In this paper we propose a unified approach for realization of mass transfer toward ends by diffeomorphisms of $M$.  
This argument together with Moser's theorem 
enables us to deduce two selection theorems for the groups ${\cal D}^+(M)$ and ${\cal D}(M; \omega)$. 
The first one is the extension of Moser's theorem to noncompact manifolds, that is, 
the existence of sections for the orbit maps under the action of ${\cal D}^+(M)$ on the space of volume forms. 
This implies that ${\cal D}(M; \omega)$ is a strong deformation retract of the group ${\cal D}^+(M; E^\omega_M)$ consisting of $h \in {\cal D}^+(M)$ which preserves the set $E^\omega_M$ of $\omega$-finite ends of $M$.   
The second one is related to the mass flow toward ends under volume-preserving diffeomorphisms of $M$.  
Let ${\cal D}_{E_M}(M; \omega)$ 
denote the subgroup consisting of all $h \in {\cal D}(M; \omega)$ which fix the ends $E_M$ of $M$. 
S.~R.~Alpern and V.\,S.\,Prasad 
introduced the topological vector space ${\cal S}(M; \omega)$ of end charges of $M$ and 
the end charge homomorphism $c^\omega : {\cal D}_{E_M}(M; \omega) \to {\cal S}(M; \omega)$, 
which measures the mass flow toward ends induced by each $h \in {\cal D}_{E_M}(M; \omega)$. 
 We show that the homomorphism $c^\omega$ has a continuous section. 
This induces the factorization ${\cal D}_{E_M}(M; \omega) \cong {\rm ker}\,c^\omega \times {\cal S}(M; \omega)$ 
and it implies that ${\rm ker}\,c^\omega$ is a strong deformation retract of ${\cal D}_{E_M}(M; \omega)$. 
\end{abstract}


\section{Introduction}

The purpose of this article is to study topological properties of 
groups of volume-preserving diffeomorphisms of noncompact manifolds. 
The group of diffeomorphisms of a manifold $M$ acts on the space of volume forms on $M$. 
We can use a family of diffeomorphisms of $M$ to transfer volumes and thereby deform a family of volume forms on $M$. 
In this article we propose a unified approach for the general problem of 
realizing data of volume transfer toward ends by diffeomorphisms of $M$ (Theorem~\ref{thm_deform}). 
We use this realization result, together with Moser's theorem, 
in order to obtain two main selection theorems related to 
volume-preserving diffeomorphisms of a noncompact manifold $M$ 
(Theorems~\ref{thm_v} and \ref{thm_a}) and 
deduce some conclusions on the group of volume-preserving diffeomorphisms of $M$ (Corollaries~\ref{cor_v} and \ref{cor_a}). 

Suppose $M$ is a connected oriented separable metrizable smooth $n$-manifold possibly with boundary and 
$\omega$ is a positive volume form on $M$. 
Let ${\cal D}^+(M)$ denote the group of orientation-preserving diffeomorphisms of $M$ endowed with the compact-open $C^\infty$-topology 
and ${\cal D}(M; \omega)$ denote the subgroup consisting of $\omega$-preserving diffeomorphisms of $M$. 

Our first concern is the relation between the groups ${\cal D}^+(M) \supset {\cal D}(M; \omega)$. 
Let ${\cal V}^+(M)$ denote the space of positive volume forms $\mu$ on $M$   
endowed with the compact-open $C^\infty$-topology, and 
for $m \in (0, \infty]$ let  
${\cal V}^+(M, m) = \{ \mu \in {\cal V}^+(M) \mid \mu(M) = m\}$.  
The group ${\cal D}^+(M)$ acts continuously 
on ${\cal V}^+(M, m)$ by $h \cdot \mu = h_\ast \mu \ (= (h^{-1})^\ast \mu)$ and 
the group ${\cal D}(M; \omega)$ coincides with the stabilizer of $\omega$ under  this action. 

When $M$ is compact, Moser's theorem \cite{Bany1, Bany2, Mos} asserts the transitivity and  
the existence of sections of the orbit maps under this action. 
Similar results for measure-preserving homeomorphisms were obtained by von Neumann-Oxtoby-Ulam \cite{OU} 
and A.\,Fathi \cite{Fa}. 

When $M$ is noncompact, it is necessary to include some informations on the ends of $M$ \cite{Be1, GS}. 
For an open subset $F$ of the space $E_M$ of the ends of $M$, 
consider the subgroup ${\cal D}^+(M; F) = \{ h \in {\cal D}^+(M) \mid \overline{h}(F) = F\}$, 
where $\overline{h}$ is the natural homeomorphic extension of $h$ to the end compactification of $M$. 
We also need to introduce the spaces 
\[ \mbox{${\cal V}^+(M; F)_{ew} = \{ \mu \in {\cal V}^+(M) \mid E_M^\mu = F \}$ \ \ and \ \ 
${\cal V}^+(M; m, F)_{ew} = {\cal V}^+(M; F) \cap {\cal V}^+(M; m)$, } \] 
where $E_M^\mu$ is the set of $\mu$-finite ends of $M$. 
According to R.\,Berlanga \cite{Be1} these spaces are endowed with the finite-ends  weak $C^\infty$-topology $ew$ (cf.\ Section 2.4). 
The group ${\cal D}^+(M; F)$ acts continuously 
on ${\cal V}^+(M; m, F)_{ew}$ by $h \cdot \mu = h_\ast \mu$, 
and the subgroup ${\cal D}(M; \omega)$ coincides with the stabilizer of $\omega$ under this action. 
The transitivity of this action was shown by R.\,E.\,Greene - K.\,Shiohama \cite{GS}. 
The similar problem for measure-preserving homeomorphisms were studied by R.\,Belranga \cite{BE}. 
He has really obtained the result on existence of sections for orbit maps \cite{Be1}. 
The first application of Theorem~\ref{thm_deform} is the $C^\infty$-version of this result, that is,  
the extension of Moser's theorem to noncompact manifolds. 
Let ${\cal D}_\partial(M)_1$ denote the path-component of $id_M$ in the group ${\cal D}_\partial(M) = \{ h \in {\cal D}(M) \mid h|_{\partial M} = id_{\partial M} \}$. 

\begin{theorem}\label{thm_v}
Suppose $P$ is any topological space and 
$\mu, \nu : P \to {\cal V}^+(M; F)_{ew}$ are continuous maps with $\mu_p(M) = \nu_p(M)$ $(p \in P)$. 
Then there exists a continuous map $h : P \to {\cal D}_\partial(M)_1$ such that for each $p \in P$ \ \ {\rm (1)} ${h_p}_\ast \mu_p = \nu_p$ \ \ and \ \ {\rm (2)} if $\mu_p = \nu_p$, then $h_p = id_M$. 
\end{theorem}

\begin{cor}\label{cor_v} Let $\omega \in {\cal V}^+(M)$ and set $m = \omega(M)$ and $F = E_M^\omega$. 
\begin{itemize} 
\item[(1)] The orbit map 
$\pi_\omega : {\cal D}^+(M;F) \to {\cal V}^+(M; m, F)_{ew}$, 
$\pi_\omega(h) = h_\ast \omega$, 
admits a continuous section 
$\sigma : {\cal V}^+(M; m, F)_{ew} \to {\cal D}_\partial(M)_1$ such that 
$\sigma(\omega) = id_M$. 
\item[(2)] \,{\rm (i)} ${\cal D}^+(M;F) \cong {\cal V}^+(M; m, F)_{ew} \times {\cal D}(M; \omega)$.\\ 
{\rm (ii)} ${\cal D}(M; \omega)$ is a strong deformation retract of ${\cal D}^+(M;F)$.   
\end{itemize} 
\end{cor} 

Next we are concerned with an internal structure of the group ${\cal D}(M; \omega)$. 
When $M$ is noncompact, one can measure 
volume transfer toward ends of $M$ under volume-preserving diffeomorphisms which fix the ends $E_M$.  
This quantity is described by the end charge homomorphism 
$$c^\omega : {\cal D}_{E_M}(M; \omega) \to {\cal S}(M; \omega)$$ 
introduced by S.~R.~Alpern and V.\,S.\,Prasad \cite{AP}. 
An end charge of $M$ is a finitely additive signed measure 
on the algebra of clopen subsets of $E_M$. 
Let ${\cal S}(M)$ denote the topological linear space of all end charges of $M$ with the weak topology and 
let ${\cal S}(M; \omega)$ denote the linear subspace of ${\cal S}(M)$ consisting of end charges $c$ of $M$ with $c(E_M) = 0$ and $c |_{E_M^\omega} = 0$. 
For each $h \in {\cal D}_{E_M}(M; \omega)$ an end charge $c^\omega_h \in {\cal S}(M; \omega)$ is defined by 
\[ c^\omega_h(E_C) = \omega(C - h(C)) - \omega(h(C) - C), \]
where $C$ is any $n$-submanifold of $M$ such that ${\rm Fr}_M C$ is compact 
and $E_C \subset E_M$ is the set of ends of $C$. 
This quantity represents the total $\omega$\,-\,volume (or mass) transfered by $h$ into $C$ and into $E_C$ in the last. 
Hence, the end charge $c^\omega_h$ describes mass flow toward ends induced by $h$. 

The group ${\cal D}_{E_M}(M; \omega)$ acts continuously on 
${\cal S}(M; \omega)$ by 
$h \cdot a = c^\omega_h + a$. 
The end charge homomorphism $c^\omega : {\cal D}_{E_M}(M; \omega) \to S(M; \omega)$ coincides with the orbit map at $0 \in {\cal S}(M; \omega)$ under this action. 
In \cite{Ya} we have shown that the end charge homomorphism for measure-preserving homeomorphisms has a continuous section. 
As the second application of Theorem~\ref{thm_deform}, we show that 
the end charge homomorphism $c^\omega$ for 
volume-preserving diffeomorphisms also has a continuous (non-homomorphic) section. 

\begin{theorem}\label{thm_a} 
Suppose $P$ is any topological space and 
$\mu : P \to {\cal V}^+(M)$ and $a : P \to {\cal S}(M)$ are continuous maps such that $a_p \in {\cal S}(M; \mu_p)$ $(p \in P)$. 
Then there exists a continuous map $h : P \to {\cal D}_\partial(M)_1$ 
such that for each $p \in P$ \\[1mm]  
\hspace{20mm} {\rm (1)} $h_p \in {\cal D}_\partial(M; \mu_p)_1$, \ 
{\rm (2)} $c_{h_p}^{\mu_p} = a_p$ \ and \ 
{\rm (3)} if $a_p = 0$, then $h_p = id_M$.  
\end{theorem}

\begin{cor}\label{cor_a} Let $\omega \in {\cal V}^+(M)$.  
\mbox{} 
\begin{itemize}
\item[(1)] The end charge homomorphism 
$c^\omega : {\cal D}_{E_M}(M; \omega) \to {\cal S}(M; \omega)$ 
admits a continuous section 
$s : {\cal S}(M; \omega) \to {\cal D}_\partial(M; \omega)_1$ 
such that 
$s(0) = id_M$. 
\item[(2)] 
\begin{itemize}
\item[(i)\,] ${\cal D}_{E_M}(M; \omega) 
\cong {\rm ker}\,c^\omega \times {\cal S}(M; \omega)$. 
\item[(ii)] ${\rm ker}\,c^\omega$ is a strong deformation retract of 
${\cal D}_{E_M}(M; \omega)$. 
\end{itemize}
\end{itemize}
\end{cor}

The group ${\rm ker}\,c^\omega$ contains the subgroup
${\cal D}^c(M; \omega)$ consisting of 
$\omega$-preserving diffeomorphisms with compact support. 
In a succeeding paper we will investigate 
the relation between these subgroups.  

This paper is organized as follows.  
Section 2 is devoted to generalities on 
end compactifications, diffeomorphism groups, spaces of volume forms and 
end charge homomorphism. 
In Sections 3.1 and 3.2 we discuss volume transfer by engulfing isotopies and obtain a fundamental deformation lemma. In Section 3.4 we state and prove Theorem~\ref{thm_deform}, the main theorem in this article. 
Theorems~\ref{thm_v}, ~\ref{thm_a} and their corollaries are proved in Section 4.  



\section{Preliminaries} 

\subsection{Conventions} \mbox{} 

Throughout the paper, a $C^\infty$ $n$-manifold $M$ means 
a separable metrizable $C^\infty$ $n$-manifold possibly with boundary. 
We add the phrase ``possibly with corners'' 
when we allow that $M$ has corners modeled on $[0, \infty)^2 \times {\Bbb R}^{n-2}$. 
In this paper corners appear only when we cut an $n$-manifold with boundary 
by a proper $(n-1)$-submanifold. 

Suppose $M$ is a $C^\infty$ $n$-manifold possibly with corners. 
The notations ${\rm Int}\,M$  and $\partial = \partial M$ 
denote the interior and boundary of $M$ as a manifold. 
For a subset $A$ of $M$, let 
${\rm Fr}_M A$, ${\rm cl}_M A$ and ${\rm Int}_M A$ 
denote the frontier, closure and interior of $A$ relative to $M$. 

A $C^\infty$ $n$-{\em submanifold} $N$ of $M$ means a closed subset of $M$ 
such that (a) ${\rm Fr}_M N$ is a proper $C^\infty$ $(n-1)$-submanifold of $M$ and (b) ${\rm Fr}_M N$ does not meet the corners of $M$. 
Hence $N$ is a $C^\infty$ $n$-manifold with corners ${\rm Fr}_M N \cap \partial M$ (together with other corners of $M$ included in $N$) and 
$N^c \equiv cl_M(M - N)$ is also a $C^\infty$ $n$-submanifold of $M$. 

When $M$ has corners, a collar $E$ of $\partial M$ in $M$ is locally modeled on $E_0 \times {\Bbb R}^{n-2}$, where $E_0 = [0,2]^2 \setminus [0,1)^2 \subset {\Bbb R}^2$. 
The collar $E_0$ has the base $C = [0,2]^2 \setminus [0,2)^2$ and 
the collar arcs $I_x = \{ (1-t/2)x \mid t \in [0, 1] \}$ $(x \in C)$ 
(so that the level sets consist of $F_t = (1-t/2) C$ $(t \in [0,1])$).  
The notation $E = \partial M \times [0,1]$ simply indicates the (non-diffeomorphic)  parametrization of $E$. 

The symbols ${\cal B}(M)$, ${\cal K}(M)$ and ${\cal C}(M)$ denote 
the sets of Borel subsets, compact subsets and connected components of $M$ respectively. 
Let ${\cal B}_c(M) = \{ C \in {\cal B}(M) \mid {\rm Fr}_M\,C : \mbox{compact} \}$ and ${\cal B}_{c, 0}(M) =  \{ C \in{\cal B}_c(M) \mid {\rm Fr}_M\,C \ \text{has zero Lebesgue measure} \}$. 

A Radon measure on $M$ is a Borel measure $\mu$ on $M$ such that 
$\mu(K) < \infty$ for each $K \in {\cal K}(M)$. 
A Radon measure $\mu$ on $M$ is called good if $\mu(p) = 0$ for each point $p \in M$ and 
$\mu(U) > 0$ for any nonempty open subset $U$ of $M$. 
Let ${\cal M}_g^\partial(M)$ denote the space of good Radon measures $\mu$ on $M$ with $\mu(\partial M) = 0$ endowed with the weak topology. 
For $\mu, \nu \in {\cal M}_g^\partial(M)$ we say that $\nu$ is $\mu$-biregular if 
$\mu$ and $\nu$ have the same collection of null sets. 
Let ${\cal M}_g^\partial(M, \widehat{\mu}\mbox{-reg}) 
= \big\{ \nu \in {\cal M}_g^\partial(M) \mid \mbox{$\nu$ is $\mu$-biregular}\big\}$. (cf. \cite{AP, Be1, Fa, Ya}) 

By $C^\infty(M)$ we denote the space of $C^\infty$-functions $f : M \to {\Bbb R}$ endowed with the compact-open $C^\infty$-topology. For $I \subset {\Bbb R}$ we have the subspace $C^\infty(M, I) = \{ f \in C^\infty(M) \mid f(M) \subset I \}.$

The symbol $\Omega^n(M)$ denotes the space of $n$-forms on $M$ 
endowed with the compact-open $C^\infty$-topology. 
\vspace{1mm} 
When $M$ is oriented, 
${\cal V}^+(M)$ denotes the subspace of $\Omega^n(M)$ 
consisting of positive volume forms on $M$. 
For any $\omega \in {\cal V}^+(M)$, the space ${\cal V}^+(M)$ admits 
a canonical affine contraction onto $\{ \omega \}$:  
\[ \mbox{$\phi : {\cal V}^+(M) \times [0,1] \to {\cal V}^+(M)$, \ 
$\phi_t(\mu) = (1-t)\mu + t\omega$.} \]  
Each $\omega \in {\cal V}^+(M)$ induces a Radon measure $\widehat{\omega} \in  {\cal M}_g^\partial(M)$ defined by 
\[ \widehat{\omega}(B) = \int_B \,\omega \equiv \int_M \chi_B \,\omega \qquad 
(B \in {\cal B}(M)), \] 
where $\chi_B$ is the characteristic function of $B$ \   
(i.e., $\chi_B \equiv 1$ on $B$ and $\chi_B \equiv 0$ on $M - B$). 
Since $\widehat{\mu} \in {\cal M}_g^\partial(M, \widehat{\omega}\mbox{-reg})$ for any $\mu \in {\cal V}^+(M)$, 
this determines a canonical continuous map 
${\cal V}^+(M) \to {\cal M}_g^\partial(M, \omega\mbox{-reg})$. 
For simplicity, we denote $\widehat{\omega}$ by the same simbol $\omega$. 
Note that $C \in {\cal B}(M)$ has zero Lebesgue measure iff $\omega(C) = 0$. 

By ${\cal D}(M)$ we denote the group of diffeomorphisms of $M$ onto itself endowed with the compact-open $C^\infty$-topology. 
The support of $h \in {\cal D}(M)$ is defined by \ ${\rm supp}\,h = {\rm cl}_M\,\{ x \in M \mid h(x) \neq x \}$.
Let ${\cal D}^c(M) = \{ h \in {\cal D}(M) \mid {\rm supp}\,h :\,\mbox{compact}\}$ and for $A, B \subset M$ we set 
\[ \mbox{${\cal D}_A(M, B) = \big\{ h \in {\cal D}(M) \mid h|_A = id_A, \,h(B) = B \big\}$ \ \ 
and \ \ ${\cal D}_A^c(M, B) = {\cal D}_A(M, B) \cap {\cal D}^c(M)$.} \] 
When $M$ is oriented, 
${\cal D}^+(M)$ denotes the subgroup of ${\cal D}(M)$ 
consisting of orientation-preserving diffeomorphisms of $M$. 
For $\omega \in {\cal V}^+(M)$ 
we say that $h \in {\cal D}(M)$ preserves $\omega$ if the pullback $h^\ast \omega = \omega$ 
(or the pushforward $h_\ast \omega \equiv (h^{-1})^\ast \omega = \omega$). 
Let ${\cal D}(M; \omega)$ denote the subgroup of ${\cal D}^+(M)$ consisting of $\omega$-preserving diffeomorphisms of $M$.

For any subgroup $G$ of ${\cal D}(M)$,  
the symbol $G_1$ denotes the path-component of $id_M$ in $G$. 
When $G \subset {\cal D}^c(M)$, by $G_1^\ast$ we denote the subgroup of $G_1$ 
consisting of $h \in G$ which admits a path $h_t \in G$ ($t \in [0, 1]$) 
such that $h_0 = id_M$, $h_1 = h$ and there exists $K \in {\cal K}(M)$ with ${\rm supp}\,h_t \subset K$ ($t \in [0,1]$).   
Since $G$ is a topological group, for any continuous maps $\phi, \psi : P \to G$, 
the composition $\psi \phi : P \to G$, $(\psi \phi)_p = \psi_p \phi_p$ and 
the inverse $\phi^{-1} : P \to G$, $(\phi^{-1})_p = \phi^{-1}_p$ are also continuous maps. 
We say that a map $\phi : P \to G$ (or a map $\mu : P \to \Omega^n(M)$) has locally common compact support in $V \subset M$
if for each point $p \in P$ there exists a neighborhood $U$ of $p$ in $P$ and $K \in {\cal K}(V)$ such that ${\rm supp}\,\phi_q \subset K$ (or ${\rm supp}\,\mu_q \subset K$) ($q \in U$). 

\subsection{Ends of manifolds}  (cf.\,\cite{Be1})

Suppose $M$ is a noncompact, connected $C^\infty$ $n$-manifold possibly with corners. 
An end of $M$ is a function $e$ which assigns an $e(K) \in {\cal C}(M - K)$ to each $K \in {\cal K}(M)$ such that 
$e(K_1) \supset e(K_2)$ if $K_1 \subset K_2$. 
The set of ends of $M$ is denoted by $E_M$. 
The end compactification of $M$ is the space $\overline{M} = M \cup E_M$ 
endowed with the topology defined by the following conditions: 
(a) $M$ is an open subspace of $\overline{M}$, 
(b) the fundamental open neighborhoods of $e \in E_M$ are given by
\[ N(e, K) = e(K) \,\cup \,\{ e' \in E_M \mid e'(K) = e(K)\} \hspace{5mm}  
(K \in {\cal K}(M)). \]
Then, (i) $\overline{M}$ is a connected, locally connected, compact, metrizable space, (ii) $M$ is a dense open subset of $\overline{M}$ and 
(iii) $E_M$ is a compact 0-dimensional subset of $\overline{M}$. 
Thus, for any metric $d$ on $\overline{M}$ and 
any $\varepsilon > 0$ there exists a neighborhood $U$ of $E_M$ in $\overline{M}$ such that 
${\rm diam}_d\, C  < \varepsilon$ $(C \in {\cal C}(U))$. 

Let ${\rm Map}(\overline{M})$ denote 
the space of continuous maps $f : \overline{M} \to \overline{M}$ and  
let ${\cal H}(\overline{M})$ denote the group of homeomorphisms $h$ of $\overline{M}$ onto itself. 
These spaces are endowed with the the compact-open topology. 
Any metric $d$ on $\overline{M}$ induces the sup-metric $d$ on ${\rm Map}(\overline{M})$, 
which is compatible with the compact-open topology. 
For $A \subset \overline{M}$ we set 
\[ \mbox{${\rm Map}_A(\overline{M}) = \{ f \in {\rm Map}(\overline{M}) \mid f|_A = id_A \}$ \ \ and \ \ 
${\cal H}_A(\overline{M}) = \{ h \in {\cal H}(\overline{M}) \mid h|_A = id_A \}$.} \] 

For $h \in {\cal D}(M)$ and $e \in E_M$ we define $h(e) \in E_M$ by  
\[ \mbox{$h(e)(K) = h(e(h^{-1}(K))) \in {\cal C}(M - K)$ \ ($K \in {\cal K}(M)$).} \] 
Then $h$ has a unique extension $\overline{h} \in {\cal H}(\overline{M})$ defined 
by $\overline{h}(e) = h(e)$ ($e \in E_M$) and 
the map \ ${\cal D}(M) \to {\cal H}(\overline{M})$ : $h \mapsto \overline{h}$ \ is a continuous group homomorphism. 
If $h \in {\cal D}(M)_1$, then $\overline{h} \in {\cal H}_{E_M}(\overline{M})$. 
For $A, B \subset \overline{M}$ we set 
${\cal D}_A(M, B) = \big\{ h \in {\cal D}(M) \mid \overline{h}|_A = id_A, \overline{h}(B) = B\big\}$. 

For each $C \in {\cal B}_c(M)$ we set 
\[ \overline{C} = C \cup E_C, \hspace{10mm}  
E_C = \{ e \in E_M \mid e(K) \subset C \mbox{ for some } K \in {\cal K}(M)\}. \]  
Then, $E_C$ is open and closed in $E_M$ and $\overline{C}$ is a neighborhood of $E_C$ in $\overline{M}$.   

For $A, B \in {\cal B}(M)$ we write $A \sim_c B$ 
if the symmetric difference $A \Delta B = (A - B) \cup (B - A)$ is relatively compact (i.e., has the compact closure) in $M$. This is an equivalence relation and 
for $C, D \in {\cal B}_c(M)$ we have (i) $C \sim_c D$ iff $E_C = E_D$ and 
(ii) $C \sim_c h(C)$ for any $h \in {\cal D}_{E_M}(M)$. 

For $\omega \in {\cal M}_g^\partial(M)$ (or $\omega \in {\cal V}^+(M)$ when $M$ is oriented),  
we say that an end $e \in E_M$ is $\omega$-finite if $\omega(e(K)) < \infty$ for some $K \in {\cal K}(M)$. 
We set $E_M^\omega = \{ e \in E_M \mid e : \mbox{$\omega$-finite}\}$. 
Then $E_M^\omega$ is an open subset of $E_M$ and for $C \in {\cal B}_c(M)$ 
we have that $\omega(C) < \infty$ iff $E_C \subset E_M^\omega$. 

\subsection{Volume transfer and end charge homomorphism} \mbox{}  

First we introduce a quantity which measures volume transfer by diffeomorphisms (cf.\,\cite[Section 3.2]{Ya}). 
Suppose $M$ is a connected oriented $C^\infty$ $n$-manifold possibly with corners and 
$\omega \in {\cal V}^+(M)$. 
For $A, B \in {\cal B}(M)$ we write $A \sim_\omega B$ if $\omega(A \Delta B) < \infty$. 
This is an equivalence relation and $A  \sim_c B$ implies $A \sim_\omega B$. If $A \sim_\omega B$, then  
we can define the following quantity: 
\[ J^\omega(A, B) = \omega(A - B) - \omega(B - A) \ \in \ {\Bbb R}. \]
If $C \in {\cal B}_c(M)$ and $h \in {\cal D}_{E_M}(M)$, then 
$J^\omega(h^{-1}(C), C)$ is just the total $\omega$\,-\,mass transfered  into $C$ by $h$. If $A \in {\cal B}(M)$ and 
$f, g \in {\cal D}^c(M)$,  then \ \ 
$J^\omega(f^{-1}(A), g^{-1}(A)) = (f_\ast \omega - g_\ast \omega)(A)$. 

The quantity $J^\omega$ has the following formal properties: 

\begin{lemma} \label{lem-difference of volume} 
Suppose $A, B, C, D \in {\cal B}(M)$. 
\begin{itemize}
\item[(1)] If $A \sim_\omega B$, then {\rm (i)} $\omega(A) = J^\omega(A, B) + \omega(B)$ and \\ 
{\rm (ii)} if $\omega(A) < \infty$, then $\omega(B) < \infty$ and $J^\omega(A, B) = \omega(A) - \omega(B)$. 

\item[(2)] If $A \sim_\omega B \sim_\omega C$, then 
$J^\omega(A, B) + J^\omega(B, C) = J^\omega(A, C)$. 

\item[(3)] If $A \sim_\omega C$, $B \sim_\omega D$, then {\rm (i)} $A \cup B \sim_\omega C \cup D$ and \\ 
{\rm (ii)} if $A \cap B  = C \cap D = \emptyset$, then $J^\omega(A \cup B, C \cup D) = J^\omega(A, C) + J^\omega(B, D)$. 

\item[(4)] If $h \in {\cal D}(M)$ and $A \sim_{h_\ast \omega} B$, then 
$h^{-1}(A) \sim_\omega h^{-1}(B)$ and 
$J^{h_\ast \omega}(A, B) = J^\omega(h^{-1}(A), h^{-1}(B))$. 
\end{itemize}
\end{lemma} 

\begin{lemma}\label{lem_conti} Suppose $A, B \in {\cal B}_{c, 0}(M)$ and $A \sim_c B$.  
Then the following function is continuous : 
$$\Phi : {\cal V}^+(M) \times {\cal D}_{E_M}(M)^2 \lra {\Bbb R} : \ \ \Phi(\mu, f, g) =  J^\mu(f(A), g(B)).$$
\end{lemma} 

Lemma~\ref{lem-difference of volume} is a special case of \cite[Lemma 3.1]{Ya}, while Lemma~\ref{lem_conti} follows directly from \cite[Lemma 3.2]{Ya} since 
the canonical maps 
${\cal V}^+(M) \to {\cal M}_g^\partial(M, \omega\mbox{-reg})$ and 
${\cal D}_{E_M}(M) \to {\cal H}_{E_M}(M, \omega\mbox{-reg})$ are continuous. 

Next we recall basic properties of the end charge homomorphism defined in \cite[Section 14]{AP}. 
Suppose $M$ is a noncompact connected oriented $C^\infty$ $n$-manifold possibly with corners 
and $\omega \in {\cal V}^+(M)$. 
The symbol ${\cal Q}(E_M)$ denotes the algebra of clopen subsets of $E_M$. 
An {\em end charge} of $M$ is a finitely additive signed measure $c$ on ${\cal Q}(E_M)$, that is, 
a function $c : {\cal Q}(E_M) \to {\Bbb R}$ which satisfies the following condition:  
\[ \mbox{$c(F \cup G) = c(F) + c(G)$ \ for \ $F, G \in {\cal Q}(E_M)$ \ with \ $F \cap G = \emptyset$.} \]  
Let ${\cal S}(M)$ denote the space of end 
charges $c$ of $M$ with the {\em weak topology} (or the product topology). 
This topology is the weakest topology such that the function \ \ 
${\cal S}(M) \to {\Bbb R}  : \ c \mapsto c(F)$ \ \ 
is continuous for any $F \in {\cal Q}(E_M)$. 
We set 
\[ {\cal S}(M; \omega) = 
\big\{ c \in {\cal S}(M) \mid \text{\rm (i) } c(E_M) = 0, \ 
\text{(ii) } c(F) = 0 \ (F \in {\cal Q}(E_M), F \subset E_M^\omega) \}. \]
Then ${\cal S}(M)$ is a topological linear space and ${\cal S}(M; \omega)$ is a linear subspace of ${\cal S}(M)$. 

For $h \in {\cal D}_{E_M}(M; \omega)$ 
the end charge $c_h^\omega \in {\cal S}(M; \omega)$ is defined as follows: 
For any $F \in {\cal Q}(E_M)$ there exists $C \in {\cal B}_c(M)$ with $E_C = F$.  
Since $C \sim_c h(C)$, we can define as 
\[ c_h^\omega(F) = J^\omega(h^{-1}(C), C) \ 
\big(= J^\omega(C, h(C)) = \omega(C - h(C)) - \omega(h(C) - C)\big) \in {\Bbb R}.  \] 
This quantity is independent of the choice of $C$. 

\begin{prop} The map \ $c^\omega : {\cal D}_{E_M}(M; \omega) \lra {\cal S}(M; \omega)$ \ is a continuous group homomorphism \\ 
$($\cite[Section 14.9, Lemma 14.21\,(iv)]{AP}$)$. 
\end{prop} 

\subsection{Action of diffeomorphism groups on spaces of volume forms} \mbox{} 

Suppose $M$ is an oriented $C^\infty$ $n$-manifold possibly with corners. 
The group ${\cal D}^+(M)$ acts continuously 
on ${\cal V}^+(M)$ by $h \cdot \mu = h_\ast \mu$. 
The subgroup ${\cal D}(M; \omega)$ coincides with the stabilizer of $\omega \in {\cal V}^+(M)$.
For $m \in (0, \infty]$ 
we obtain the invariant subspace ${\cal V}^+(M; m) = \{ \mu \in {\cal V}^+(M) \mid \mu(M) = m \}$. For $\omega \in {\cal V}^+(M; m)$ we have 
the orbit map $\pi_\omega : {\cal D}^+(M) \to {\cal V}^+(M; m)$, $\pi_\omega(h) = h_\ast \omega$. 
Moser's theorem asserts that this orbit map has a continuous section 
in case $M$ is a compact connected oriented $C^\infty$ $n$-manifold (J.\,Moser \cite{Mos}, A.\,Banyaga \cite{Bany1, Bany2}). 
In this article we apply the next version. 
The following notations are used for a collar $E = S \times [a,b]$ and $A \subset S$, \\ 
\qquad $\begin{array}[t]{l}
E_A = \{ (x,t) \in E \mid x \in A \}, \quad E_A^+ = \{ (x,t) \in E_A \mid t \geq 0 \},  
\quad E_A^- = \{ (x,t) \in E_A \mid t \leq 0 \} \quad \text{and} \\[1.5mm] 
E_A [\delta, \e] = \{ (x,t) \in E_A \mid t \in [\delta(x), \e(x)] \} \quad 
\text{for functions $\delta, \e : S \to [a,b]$.}  
\end{array}$
\vskip 1.5mm 

\begin{theorem}\label{thm_Moser} 
Suppose $M$ is a compact connected oriented $C^\infty$ $n$-manifold possibly with corners and $E = \partial M \times [0,1]$ is a collar of $\partial M$. 
Suppose $\mu, \nu : P \to {\cal V}^+(M)$ are continuous maps and 
$\varepsilon : P \to (0, 1/2)$ is a continuous function such that for each $p \in P$ 
\[ \mbox{{\rm (i)} \ $\mu_p(M) = \nu_p(M)$ \ \ and \ \ 
{\rm (ii)} \ $\mu_p = \nu_p$ on $E[0, 2\varepsilon_p]$.} \]  
Then there exists a continuous map $\phi : P \to {\cal D}_\partial(M)_1$ 
such that for each $p \in P$ 
\[ \mbox{\rm (1) ${\phi_p}_\ast \mu_p = \nu_p$, \ \ 
(2) $\phi_p = id_M$ on $E[0, \varepsilon_p]$, \ \ 
(3) if $\mu_p = \nu_p$, then $\phi_p = id_M$.} \] 
\end{theorem}

\noindent Theorem~\ref{thm_Moser} is complemented by the next lemma, 
which is a parametrized version of \cite[Lemma A2]{Mc}. 

\begin{lemma}\label{lem_collar} Let $M$ be an oriented $C^\infty$ $n$-manifold possibly with  boundary. 
\begin{itemize}
\item[(1)] 
Suppose $S$ is a proper $(n-1)$-submanifold of $M$, 
$K$ is a closed subset of $S$, $U$ is a neighborhood of $K$ in $S$ and 
$E = S \times [-1, 1]$ is a bicollar neighborhood of $S$ in $M$. 
Then, for any continuous maps $\mu, \nu : P \to {\cal V}^+(M)$ there exist continuous maps $\phi : P \to {\cal D}_{S \cup (M - E_U)}(M)_1$ and 
$\varepsilon : P \to C^\infty(S, (0,1))$ 
such that for each $p \in P$ 
\begin{itemize} 
\item[(i)\,] ${\phi_p}_\ast \mu_p = \nu_p$ \ on \ $E_K[-\varepsilon_p, \varepsilon_p]$, 
\item[(ii)] for each $x \in S$, \ \ {\rm (a)} \ $\phi_p(E_x^\pm) = E_x^\pm$ \ \ and \ \ {\rm (b)} \ if $\mu_p = \nu_p$ on $E_x^\pm$, 
then $\phi_p = id$ on $E_x^\pm$, 

in particular, if $\mu_p = \nu_p$ on $E^\pm$, 
then $\phi_p = id$ on $E^\pm$. 
\end{itemize} 
\vskip 1mm

\item[(2)] Suppose $K$ is a closed subset of $\partial M$, $U$ is a neighborhood of $K$ in $\partial M$ and $E = \partial M \times [0, 1]$ is a collar of $\partial M$ in $M$. 
Then, for any continuous maps $\mu, \nu : P \to {\cal V}^+(M)$ there exist continuous maps $\psi : P \to {\cal D}_{\partial M \cup (M - E_U)}(M)_1$ and $\varepsilon : P \to C^\infty(\partial M, (0,1))$ 
such that for each $p \in P$, 
\begin{itemize}
\item[(i)\,] ${\psi_p}_\ast \mu_p = \nu_p$ \ on \ $E_K[0, \varepsilon_p]$, 
\item[(ii)] for each $x \in \partial M$, \ {\rm (a)} \ $\psi_p(E_x) = E_x$ \ \ and \ \ {\rm (b)} \ 
if $\mu_p = \nu_p$ on $E_x$, then $\psi_p = id$ on $E_x$, \\
in particular, if $\mu_p = \nu_p$ on $E$, then $\psi_p = id_M$.
\end{itemize} 

\end{itemize}
\end{lemma} 

\begin{remark} 
In Lemma~\ref{lem_collar} (1), the map $\phi_p$ does not fix $\partial M$ pointwise. 
To remedy this point, we can apply (2) before (1) to obtain the following conclusion. 
\begin{itemize} 
\item[(3)] Suppose $S$ is a proper $(n-1)$-submanifold of $M$ and 
$E = S \times [-1, 1]$ is a bicollar neighborhood of $S$ in $M$. 
Then, for any continuous maps $\mu, \nu : P \to {\cal V}^+(M)$ there exist continuous maps 
$\phi' : P \to {\cal D}_{\partial M \cup (M - E)}(M, S)_1$ and 
$\varepsilon' : P \to C^\infty(S, (0,1))$ 
such that for each $p \in P$, 
\begin{itemize}
\item[(i)\,] ${\phi'_p}_\ast \mu_p = \nu_p$ \ on \ $E[-\varepsilon'_p, \varepsilon'_p]$ \ \ and    
\item[(ii)] if $\mu_p = \nu_p$, then $\phi'_p = id_M$. 
\end{itemize} 
\end{itemize}
\end{remark} 

For completeness we include the proofs of 
Theorem~\ref{thm_Moser} and Lemma~\ref{lem_collar} in Appendix. 

To deal with the noncompact case, 
it is necessary to include the informations on the ends of manifolds (\cite{Be1, GS}). 
Suppose $M$ is a noncompact, connected oriented $C^\infty$ $n$-manifold possibly with corners. 
For $m \in (0, \infty]$ and an open subset $F$ of $E_M$, we set 
\[ {\cal V}^+(M; F) = \big\{ \mu \in {\cal V}^+(M) \mid E_M^\mu = F\big\} \hspace{3mm} \mbox{and} \hspace{3mm} 
{\cal V}^+(M; m, F) = {\cal V}^+(M; m) \cap {\cal V}^+(M; F). \] 
To distinguish topologies, by ${\cal V}^+(M; F)_w$ we denote the space ${\cal V}^+(M; F)$ endowed with the compact-open $C^\infty$-topology. 
If $C \in {\cal B}_c(M)$ and $E_C \subset F$, then we obtain the function \ \ 
$$\Phi_C  : {\cal V}^+(M; F)_w \to {\Bbb R}, \ \ \Phi_C(\mu) = \mu(C).$$ 
It is easily seen that this map is not continuous if $C$ is not compact. 
In \cite{Be1} R.\,Belranga overcame this problem by introducing a stronger topology called the finite-ends weak topology (on the spaces of Radon measures). Consider the subspace $M \cup F \subset \overline{M}$. 

The finite-ends weak topology on the space \ 
${\cal M}_g^\partial(M; F) = \{ \mu \in {\cal M}_g^\partial(M) \mid E_M^\mu = F\}$ \  
is the weakest topology such that 
\begin{itemize}
\item[(0)] the function \ $\Phi_f : {\cal M}_g^\partial(M; F) \to {\Bbb R}$ : $\ds \Phi_f(\mu) = \int_M f \,d\mu$ \ is continuous \\[1.5mm]
\hspace{30mm} for any continuous function $f : M \cup F \to {\Bbb R}$ with compact support. 
\end{itemize}
Let ${\cal M}_g^\partial(M; F)_{ew}$ denote 
the space ${\cal M}_g^\partial(M; F)$ endowed with the finite-ends weak topology. 

\begin{defi}\label{def_ew}
The finite-ends weak $C^\infty$ topology on ${\cal V}^+(M; F)$ is 
the weakest topology such that 
\begin{itemize}
\item[(i)\,] the identity map $id : {\cal V}^+(M; F) \to {\cal V}^+(M; F)_w$ is continuous,  
\vskip 1mm 
\item[(ii)] the function \ $\Phi_f : {\cal V}^+(M; F) \to {\Bbb R}$ : $\ds \Phi_f(\mu) = \int_M f \mu$ \ is continuous \\[1.5mm]
\hspace{30mm} for any continuous function $f : M \cup F \to {\Bbb R}$ with compact support. 
\end{itemize}
\end{defi}

\noindent Let ${\cal V}^+(M; F)_{ew}$ denote 
the space ${\cal V}^+(M; F)$ endowed with the finite-ends weak $C^\infty$ topology. 
From the definition it follows that the canonical map ${\cal V}^+(M; F)_{ew} \to {\cal M}_g^\partial(M; F)_{ew}$ is continuous. 

\begin{lemma}\label{lem_conti_ew} 
If $C \in {\cal B}_c(M)$ and $E_C \subset F$, 
then the function \  
$\Phi_C  : {\cal V}^+(M; F)_{ew} \to {\Bbb R}$, \ $\Phi_C(\mu) = \mu(C)$ \ is continuous. 
\end{lemma} 

\begin{proof} 
Let $K = {\rm cl}_M C \cup E_C$ and 
take an open neighborhood $U$ of ${\rm Fr}_M C$ in $M$ such that ${\rm cl}_M U$ is compact. 
Then $K, K - U \in {\cal K}(M \cup F)$ and 
$\overline{C} = C \cup E_C$ is a neighborhood of $K - U$ in $M \cup F$.
Thus, there exists a continuous function $f : M \cup F \to [0,1]$ such that 
$f = 1$ on $K - U$ and $f = 0$ on $(M \cup F) \setminus \overline{C}$. 

Since $f$ has a compact support in $K$, 
by the condition (ii) in Definition~\ref{def_ew} we have the continuous function 
$\Phi_f : {\cal V}^+(M; F)_{ew} \to {\Bbb R}$. 
Since the Borel measurable function $g = \chi_C - f : M \to [0, 1]$ has a compact support in $cl_M U$, by the condition (i) in Definition~\ref{def_ew} the map 
$\Phi_g : {\cal V}^+(M; F)_{ew} \to {\Bbb R}$ is continuous. 
Thus $\Phi_C = \Phi_{\chi_C} = \Phi_f + \Phi_g$ is also continuous. 
\end{proof}

The group ${\cal D}^+(M; F)$ acts continuously 
on ${\cal V}^+(M; F)_{ew}$ by $h \cdot \mu = h_\ast \mu$. 
The subspace ${\cal V}^+(M; m, F)_{ew}$ is invariant under this action and for $\omega \in {\cal V}^+(M; m, F)_{ew}$ 
we have the orbit map 
$\pi_\omega : {\cal D}^+(M; F) \to {\cal V}^+(M; m, F)_{ew}$, $\pi_\omega(h) = h_\ast \omega$. 
The subgroup ${\cal D}^+(M; \omega)$ coincides with the stabilizer of $\omega$. The transitivity of this action on ${\cal V}^+(M; m, F)_{ew}$ is proved in \cite{GS}. 

For the sake of notational simplicity we follow the next convention in the subsequent  sections. 
For continuous maps $\phi : P \to {\cal D}^+(M)$, 
$\mu, \nu : P \to {\cal V}^+(M)$ and an $n$-submanifold $N$ of $M$, 
we obtain continuous maps $\phi_\ast \mu : P \to {\cal V}^+(M)$, $(\phi_\ast \mu)_p = {\phi_p}_\ast \mu_p$ and 
$\nu|_N :  P \to {\cal V}^+(N)$, $(\nu|_N)_p = \nu_p|_N$ . 

\begin{notation} 
The identities for maps \ \ $(\phi_\ast\mu)(N) = \nu(N)$, \ $\phi|_{N} = id_{N}$ \  etc. mean that 
$({\phi_p}_\ast\mu_p)(N) = \nu_p(N)$ $(p \in P)$, \ $\phi_p|_{N} = id_{N}$ $(p \in P)$ \ etc. 
\end{notation}


\section{Realization theorem for volume transfer toward ends}

In this section we discuss realization problem of volume transfer toward ends by diffeomorphisms. 
Main theorem~\ref{thm_deform} is applied in the next section to deduce Theorems~\ref{thm_v} and \ref{thm_a}. 

\subsection{Volume transfer by engulfing isotopy} \mbox{} 

Suppose $M$ is a connected oriented $C^\infty$ $n$-manifold possibly with corners 
and $d$ is any metric on $\overline{M}$. 
Consider a decomposition 
$M = L \cup_S N$ such that 
\begin{itemize}
\item[(i)\,] $L$ and $N$ are connected $C^\infty$ $n$-submanifolds of $M$,   
\item[(ii)] $S \equiv L \cap N = {\rm Fr}_M L  = {\rm Fr}_M N$, which is a compact proper $C^\infty$ $(n-1)$-submanifold of $M$. 
\end{itemize} 

\begin{lemma}\label{lem_eng}
There exists a continuous map $f : (-\infty, \infty) \to {\cal D}^{c}_\partial(M)_1^\ast$ such that  
\begin{itemize} 
\item[(1)] $f_0 = id$, \ $f_s(L) \subsetneqq f_t(L)$ $(s < t)$, 
\item[(2)] there exists a subpolyhedron $F$ of $M$ $($with respect to a $C^\infty$-triangulation of $M$$)$ such that 
\begin{itemize} 
\item[(i)\ ] ${\rm dim}\,F = n-1$ and $\partial M \subset F$, 
\item[(ii)\,] the map $f$ has a locally common compact support in $M - F$ \\
$($\,i.e., for any $T > 0$ there exists $K \in {\cal K}(M-F)$ such that 
${\rm supp}\,f_t \subset K$ $(t \in [-T, T])$\,$)$, 
\item[(iii)] for any $K \in {\cal K}(M - F)$ there exist $-\infty < s < t < \infty$ with $K \subset f_t(L) - f_s(L)$, 
\end{itemize}
\item[(3)] $\{ f_t \}_{-\infty < t < \infty}$ is equi-continuous with respect to $d|_M$. 
\end{itemize}
\end{lemma} 

\begin{proof}
Take a $C^\infty$-triangulation $\tau$ of $M$ 
such that $L$ and $N$ are subcomplexes of $\tau$. 
For $k = 0, 1, \cdots, n$ let $\tau(k)$ denote the set of $k$-simplices of $\tau$ 
and for $t \in \tau$ let $\tau_t = \big\{ t \in \tau \mid t \subsetneqq s \big\}$. 
Consider the barycentric subdivision $\tau'$ of $\tau$.  
The barycenter of a simplex $s \in \tau$ is denoted by $b(s)$. 
If $s_1, s_2 \in \tau$ and $s_1 \subsetneqq s_2$,  
then we have the associated 1-simplex $\langle b(s_1), b(s_2) \rangle$ of $\tau'$. 
The dual 0-cell $e(s)$ associated with an $s \in \tau(n)$ is the barycenter $b(s)$ of $s$. 
The dual 1-cell $e(t)$ associated with an $t \in \tau(n-1)$ is the union of 
1-simplices $\langle b(t), b(s) \rangle$ ($s \in \tau_t$) of $\tau'$. 
Note that $\#\tau_t = 1$ if $t \subset \partial M$ and $\#\tau_t = 2$ otherwise.  
The dual 1-skeleton $S_M$ of $\tau$ is the union of all dual 0-cells and dual 1-cells of $\tau$. 
A maximal tree of $S_M$ is a tree $T$ of the graph $S_M$ 
which is a union of dual 1-cells of $\tau$ and includes all dual 0-cells of $\tau$. 

Consider the triangulation $\tau|_L$ of $L$ and $\tau|_N$ of $N$. 
Let $S_L$ and $S_N$ be the dual 1-skeletons of $\tau|_L$ and $\tau|_N$ and 
take maximal trees $T_L$ of $S_L$ and $T_N$ of $S_N$ such that 
$T_L \subset {\rm Int}\,L$ and $T_N \subset {\rm Int}\,N$. 
Take any $(n-1)$-simplex $s_0 \in \tau|_S$, and 
let $T_S$ denote the dual 1-cell associated to $s_0$.
Then $T = T_L \cup T_S \cup T_N$ is a maximal tree of the dual 1-skeleton $S_M$ of $\tau$. 
Let $F = \cup \{ s \in \tau \mid s \cap T = \emptyset\}$ and $U = M - F$. 
Then $F$ is a subpolyhedron of $M$ such that 
${\rm dim}\,F = n-1$, $\partial M \subset F$, $S \cap F = S - {\rm Int}\,s_0$,  
and $U$ is an open regular neighborhood of $T$, which is an open $C^\infty$ $n$-cell. 

Consider the decomposition of the end compactification $\overline{M} = \overline{L} \cup_S \overline{N}$. 
Let $\overline{F} = F \cup E_M$, $\overline{F}_L = \overline{F} \cap \overline{L}$ and 
$\overline{F}_N = \overline{F} \cap \overline{N}$. 
Using $(U, T)$, we can find an engulfing pesudo-isotopy 
$\overline{f} : [-\infty, \infty] \to {\it Map}_{\overline{F}}(\overline{M})$ 
in $\overline{M}$ such that 
\begin{itemize} 
\item[(i)\ ] 
\begin{itemize}
\item[(a)] $\overline{f}_0 = id$ and $\overline{f}_t \in {\cal H}_{\overline{F}}(\overline{M})$ ($t \in (-\infty, \infty)$), 
\item[(b)] 
for any $T \in (0, \infty)$ there exists $K \in {\cal K}(U)$ such that ${\rm supp}\,\overline{f}_t \subset K$ \ ($t \in [-T, T]$), 
\end{itemize} 
\item[(ii)\,] 
\begin{itemize}
\item[(a)] $f_t = \overline{f}_t|_M \in {\cal D}_F^c(M)$ ($t \in (-\infty, \infty)$),    
\item[(b)] the map \ $f : (-\infty, \infty) \to {\cal D}_F^c(M)_1^\ast$ : $t \mapsto f_t$ \ is continuous, 
\item[(c)] $f$ has locally common compact support in $U$, 
\end{itemize}
\item[(iii)] 
\begin{itemize}
\item[(a)] $\overline{f}_s(\overline{L}) \subsetneqq \overline{f}_t(\overline{L})$ for $-\infty \leq s < t \leq \infty$, 
\item[(b)] for any $K \in {\cal K}(U)$ there exist $-\infty < s < t < \infty$ 
such that $K \subset \overline{f}_t(\overline{L}) - \overline{f}_s(\overline{L})$, 
\item[(c)] $\overline{f}_{-\infty}(\overline{L}) = \overline{F}_L$, 
$\overline{f}_{-\infty}(\overline{N}) = \overline{M}$, and 
$\overline{f}_{\infty}(\overline{L}) = \overline{M}$, 
$\overline{f}_{\infty}(\overline{N}) = \overline{F}_N$.
\end{itemize} 
\end{itemize} 

Since the compact family $\{ \overline{f}_t \}_{-\infty \leq t \leq \infty}$ is equi-continuous with respect to $d$, 
the family $\{ f_t \}_{-\infty < t < \infty}$ is also equi-continuous with respect to $d|_M$. 
Hence, the map $f$ given in (ii) satisfies the required conditions. 
\end{proof} 
 
Now we use the isotopy $h_t = f_t^{-1}$ to transfer volume on $M$.  
Consider the continuous maps
\begin{itemize}
\item[(1)] 
$\lambda : {\cal V}^+(M) \times {\Bbb R} \to {\Bbb R}, \hspace{3mm} 
\lambda(\mu, t) = J^\mu(h_t^{-1}(L), L) =  
\begin{cases} 
\ \ \, \mu(f_t(L) - L) & (t \geq 0) \\[0.5mm]
- \mu(L - f_t(L)) & (t \leq 0) 
\end{cases}$, 
\vskip 2mm 
\item[(2)] ${\cal W}^+(M)
\ = \ \big\{ (\mu, a) \in {\cal V}^+(M) \times {\Bbb R} \ \big| \ 
a \in (- \mu(L), \mu(N)) \big\}$ (an open subset of ${\cal V}^+(M) \times {\Bbb R}$), 
\vskip 1mm 
\item[(3)] $t : {\cal W}^+(M) \to {\Bbb R}$, \ \ $t(\mu, a) = \lambda(\mu, \ast)^{-1}(a)$ \hspace{3mm} 
(i.e., $t = t(\mu, a)$ iff $a = \lambda(\mu, t)$), 
\vskip 1mm 
\item[(4)] $H : {\cal W}^+(M) \to {\cal D}_\partial^c(M)_1^\ast$, \ 
$H_{(\mu,a)} = h_{t(\mu, a)}$ 
\end{itemize} 
\vskip 1mm 

\noindent Since the map $\lambda(\mu, \ast) : {\Bbb R} \to (-\mu(L), \mu(N))$ is a monotonically increasing homeomorphism, the map $t$ is well defined. 
Note that $h_0 = id_M$ and $t(\mu, 0) = 0$. 

\begin{lemma}\label{lem_basic} The map  $H$ has the following properties:    
\begin{itemize}
\item[(i)\ ] $H$ has locally common compact support in ${\rm Int}\,M$.  
\item[(ii)\,] $\{ H_{(\mu,a)}{}^{-1} \}_{(\mu,a)}$ is equi-continuous with respect to $d|_M$.
\item[(iii)] $J^\mu(H_{(\mu,a)}^{-1}(L), L) = a$. 
\hspace{5mm} 
{\rm (iv)} $H_{(\mu,a)} = id_M$ iff $a = 0$. 
\end{itemize} 
\end{lemma} 

\begin{remark}\label{rem_homotopy} 
The continuous map 
\[ \mbox{$\widetilde{H} : {\cal W}^+(M) \times [0,1] \to {\cal D}_\partial^c(M)_1^\ast$, \ 
$\widetilde{H}_{(\mu,a), s} = H_{(\mu, sa)}$} \] 
 provides a homotopy from $\widetilde{H}_0 \equiv id_M$ to $\widetilde{H}_1 = H$ with locally common compact support in ${\rm Int}\,M$. 
\end{remark} 

\subsection{Fundamental deformation lemmas} \mbox{} 

Suppose $M$ is a connected oriented $C^\infty$ $n$-manifold possibly with corners, $d$ is any metric on $\overline{M}$ and 
$N$ is a connected $C^\infty$ $n$-submanifold 
of $M$ such that ${\rm Fr}_M N$ is compact. 
Let ${\cal C}(N^c) = \{ A_1, \cdots, A_m \}$. 

\begin{lemma}\label{lem_deform} 
Suppose $\mu : P \to {\cal V}^+(M)$ are continuous maps and 
$a(i) : P \to {\Bbb R}$ $(i = 0, 1,\cdots, m)$ are continuous functions 
such that 
\[ \mbox{
{\rm (a)} $\sum_{i=0}^m a(i) = 0$ \ \ and \ \ 
{\rm (b)} $a(0) > - \mu(N)$, $a(i) > - \mu(A_i)$ $(i = 1,\cdots, m)$.} \] 
Then there exists a continuous map $\phi : P \to {\cal D}_\partial^c(M)_1^\ast$ such that  
\begin{itemize} 
\item[(1)] 
\begin{itemize} 
\item[(i)\,] $\phi$ has locally common compact support in ${\rm Int}\,M$ 
$($and there exists a homotopy from $\phi$ to the constant map $id_M$ with locally common compact support in ${\rm Int}\,M$$)$, 
\item[(ii)] the family $\big\{ \phi_p^{-1} \big\}_p$ is equi-continuous with respect to $d|_M$, 
\end{itemize} 

\item[(2)] {\rm (i)} $J^\mu(\phi^{-1}(N), N) = a(0)$ \ and \ {\rm (ii)} 
$J^\mu(\phi^{-1}(A_i), A_i) = a(i)$ $(i = 1, \cdots, m)$, 
\item[(3)] if $p \in P$ and $a_p(i) = 0$ $(i \in \{ 1, \cdots, m \})$, then $\phi_p = id_M$. 
\end{itemize} 
\end{lemma} 

\begin{remark}\label{rem_deform} In Lemma~\ref{lem_deform} we note the following points. 
\begin{itemize} 
\item[(A)] The condition (2)(ii) implies (2)(i). 
This follows from the condition (a) and the equation 
\[ \mbox{$J^\mu(\phi^{-1}(N), N) + \sum_i J^\mu(\phi^{-1}(A_i), A_i)
 = J^\mu(\phi^{-1}(M), M) = J^\mu(M, M) = 0.$} \] 
\item[(B)] 
If the condition (b) is replaced by \hspace{5mm} 
(b)$'$ $a(0) < \mu(N)$, $a(i) < \mu(A_i)$ $(i = 1,\cdots, m)$, \\
then the condition (2) is replaced by \\ 
\hspace{5mm} (2)$'$ (i) $J^\mu(N, \phi^{-1}(N)) = a(0)$ \ and \ (ii) 
$J^\mu(A_i, \phi^{-1}(A_i)) = a(i)$ $(i = 1, \cdots, m)$. \\
This is verified by applying Lemma~\ref{lem_deform} to the functions $-a_i$. 
\end{itemize} 
\end{remark}

\begin{proof}[\bf Proof of Lemma~\ref{lem_deform}]  
{\bf [1]} First we verify the case $m = 1$. 
The argument in Section 3.1 can be applied to $M = A_1 \cup N$ so to yield the associated  map $H : {\cal W}^+(M) \to {\cal D}_\partial^c(M)_1^\ast$. 
From the conditions (a) and (b) it follows that $(\mu_p, a_p(1)) \in {\cal W}^+(M)$ (i.e., $a_p(1) \in (-\mu_p(A_1), \mu_p(N))$) for each $p \in P$. Thus, we obtain the composition 
\[ \mbox{$\phi : P \to {\cal D}_\partial^c(M)_1^\ast$, \ \ 
$\phi_p = H_{(\mu_p, a_p(1))}$ $(p \in P)$.} \]
The required properties of $\phi$ are verified by Lemma~\ref{lem_basic} and 
Remarks ~\ref{rem_homotopy} and ~\ref{rem_deform}\,(A). 

{\bf [2]} We proceed by the induction on $m$. For $m = 0$ we can take $\phi \equiv id_M$. 
Suppose $m \geq 1$ and assume that the assertion holds for $m-1$.  
To prove the assertion for $m$, first 
we apply the inductive hypothesis to the data : 
$M = (N \cup A_1) \cup A_2 \cup \cdots \cup A_m$ and $a(0) + a(1), a(2), \cdots, a(m)$ 
to yield a map $\psi : P \to {\cal D}_\partial^c(M)_1^\ast$ which satisfies the conditions (1), (3) (for $a(0) + a(1), a(2), \cdots, a(m)$) and 
\begin{itemize} 
\item[(2)$_1$] (i) $J^\mu(\psi^{-1}(N \cup A_1), N \cup A_1) = a(0) + a(1)$ \ and  \ 
(ii) $J^\mu(\psi^{-1}(A_i), A_i) = a(i)$ $(i = 2, \cdots, m)$. 
\end{itemize}
Next we apply [1] (the case $m=1$) to the data : \ 
$M_1 = N \cup A_1$, \ $(\psi_\ast \mu)|_{M_1} : P \to {\cal V}^+(M_1)$ \ \ 
and 
\[ \mbox{$b(0), b(1) : P \to {\Bbb R}$, \ \ $b(0) = a(0) - J^\mu(\psi^{-1}(N), N)$, \ \ 
$b(1) = a(1) - J^\mu(\psi^{-1}(A_1), A_1)$.} \] 
Note that \ (a) $b(0) + b(1) 
= a(0) + a(1) - J^\mu(\psi^{-1}(N \cup A_1), N \cup A_1) = 0$,  
\begin{itemize}
\item[(b)] 
$b(0) = a(0) - J^\mu(\psi^{-1}(N), N) 
> -\mu(N) - J^\mu(\psi^{-1}(N), N) = - \mu(\psi^{-1}(N)) = - (\psi_\ast \mu)(N)$, \\   
and similarly, \ $b(1) > - (\psi_\ast \mu)(A_1)$. 
\end{itemize}
Then we obtain a map $\chi : P \to {\cal D}_\partial^c(M_1)_1^\ast$ which satisfies the conditions (1), (2) and (3) for these data. 
The condition (1)\,(i) assures that this map has a canonical extension by $id$, 
$$\chi : P \to {\cal D}_{\partial \cup A_2 \cup \cdots \cup A_m}^c(M)_1^\ast,$$ 
which also satisfies the conditions (1), (3) (for $b(0), b(1)$) and 
\begin{itemize}
\item[(2)$_2$] (i) $J^{\psi_\ast \mu}(\chi^{-1}(N), N) = b(0)$ \ and  \ 
(ii) $J^{\psi_\ast \mu}(\chi^{-1}(A_1), A_1) = b(1)$. 
\end{itemize}
Finally, we see that the map 
\hspace{2mm} $\phi : P \to {\cal D}_\partial^c(M)_1^\ast$, \ \ 
$\phi_p = \chi_p \psi_p$ \hspace{2mm} satisfies the required conditions. 
\begin{itemize}
\item[(1)] Since both $\psi$ and $\chi$ satisfy the condition (1), so is the map $\phi$. 
\item[(2)] From the condition (2)$_2$\,(ii) it follows that 
\begin{align*}
a(1) &= J^\mu(\psi^{-1}(A_1), A_1) + J^{\psi_\ast \mu}(\chi^{-1}(A_1), A_1) \\[1mm] 
&= J^\mu(\psi^{-1}(A_1), A_1) + J^\mu(\psi^{-1}\chi^{-1}(A_1), \psi^{-1}A_1)
= J^\mu(\phi^{-1}(A_1), A_1). \\
\intertext{For $i = 2, \cdots, m$, since $\chi = id$ on $A_i$, 
the condition (2)$_1$\,(ii) implies that} 
a(i) & = J^\mu(\psi^{-1}(A_i), A_i) = J^\mu(\phi^{-1}(A_i), A_i).
\end{align*}
The condition (2)(i) follows from (2)(ii) and Remark~\ref{rem_deform}\,(A). 
\item[(3)] If $a_p(i) = 0$ $(i \in \{ 1, \cdots, m \})$, then we have 
(i) $\psi_p = id_M$ and (ii) $b_p(0) = b_p(1) = 0$ so that $\chi_p = id_M$. 
Thus $\phi_p = id_M$. 
\end{itemize}

\noindent This completes the proof. 
\end{proof} 

\subsection{Admissible sequences of diffeomorphisms}\mbox{} 

Suppose $M$ is a noncompact connected $C^\infty$ $n$-manifold 
possibly with boundary 
and $d$ is any fixed metric on $\overline{M}$. 
Let ${\cal N}(M)$ denote the set of compact, connected $C^\infty$ $n$-submanifold $N$ of $M$ such that each $C \in {\cal C}(N^c)$ is noncompact.  
(We assume that $\emptyset \in {\cal N}(M)$.) 
We also let   
\[ \mbox{${\cal N}^{(2)}(M) = \big\{ (K, L) \in {\cal N}(M)^2 \mid 
\text{(i) $K \subset {\rm Int}_M L$, (ii) 
$L \cap A$ is connected for each $A \in {\cal C}(K^c)$}
\big\}.$} \] 
Note that (i) for any $N \in {\cal N}(M)$ the complementary region $N^c$ 
is a $C^\infty$ $n$-submanifold of $M$ 
consisting of finitely many noncompact connected components and 
(ii) since ${\rm dim}\,E_M = 0$, for any $K \in {\cal N}(M)$ and any $\varepsilon > 0$ there exists an $L \in {\cal N}(M)$ such that $(K, L) \in {\cal N}^{(2)}(M)$ and 
${\rm diam}_{\,d}\,B < \varepsilon$ ($B \in {\cal C}(L^c)$). 

Let $P$ be a fixed topological space. 

\begin{defi}\label{def_adm_seq} A sequence $(K_k, L_k, f^k, g^k)$ $(k=1,2,\cdots)$ is said to be admissible if it satisfies the following conditions for each $k=1,2,\cdots$: \hspace{10mm} (Let $L_0 = \emptyset$ and $f^0 = g^0 \equiv id_M$)
\begin{itemize} 
\item[$(1)_k$] $K_k, L_k \in {\cal N}(M)$ and 
$(L_{k-1}, K_k), (K_k, L_k) \in {\cal N}^{(2)}(M)$.   
\vskip 1mm 
\item[$(2)_k$] 
\begin{itemize} 
\item[(i)\ ] $f^k, g^k : P \to {\cal D}_\partial^c(M)_1^\ast$ are continuous maps. 
\item[(ii)\,] $f^k = \phi^k f^{k-1}$ and $g^k = \psi^k g^{k-1}$ \\
\hspace{10mm} 
for some continuous maps $\phi^k :P \to {\cal D}_{\partial \cup L_{k-1}}^c(M)_1^\ast$ and 
$\psi^k :P \to {\cal D}_{\partial \cup K_{k}}^c(M)_1^\ast$. 
\end{itemize} 
\vskip 1mm 
\item[$(3)_k$] 
\begin{itemize} 
\item[(i)\,] $f^k$ and $g^k$ have locally common compact support in ${\rm Int}\,M$.
\item[(ii)] $\{ (f^k_p)^{-1}\}_p$ and $\{ (g^k_p)^{-1}\}_p$ are equicontinuous with respect to $d|_M$.
\end{itemize} 
\vskip 1mm 
\item[$(4)_k$] 
\hspace{-1.7mm} 
\begin{tabular}[t]{clll}
(i) & ${\rm diam}\,A \leq 2^{-k}$, & ${\rm diam}\,(g^{k-1})^{-1}(A) \leq 2^{-k}$ & $(A \in {\cal C}(K_k^c))$. \\[1mm] 
(ii) & ${\rm diam}\,B \leq 2^{-k}$, & ${\rm diam}\,(f^k)^{-1}(B) \leq 2^{-k}$ & $(B \in {\cal C}(L_k^c))$. 
\end{tabular}
\end{itemize} 
\end{defi} 

\begin{lemma}\label{lem_lim} Suppose $(K_k, L_k, f^k, g^k)_k$ is an admissible sequence. 
\begin{itemize}
\item[(1)] The sequences $(f^k)_k$ and $(g^k)_k$ converge  
$d|_M$-uniformly to some continuous maps $f, g : P \to {\cal D}_\partial(M)_1$ respectively. 
\vskip 1mm 
\item[(2)] $f^{-1}|_{L^k} = (f^{k})^{-1}|_{L^k}$, \ 
$g^{-1}|_{K^k} = (g^{k-1})^{-1}|_{K^k}$ \ $(k \geq 1)$.
\end{itemize} 
\end{lemma}

\begin{proof} We verify the statements for the map $f$.  
The argument for the map $g$ is completely same.  
Note that the metric $d$ on $\overline{M}$ induces the sup-metric $d$ on ${\rm Map}(\overline{M})$ (and on ${\rm Map}(M)$). 

[1] The sequence $(f^k)_k$ has the following properties: 
\begin{itemize}
\item[(i)\,] $\ds d(f^{k}, f^{k+1}) \leq 2^{-k}$, \ 
$\ds d\big((f^{k})^{-1}, \big(f^{k+1})^{-1}\big) \leq 2^{-k}$
 \ $(k \geq 1)$ \ \ 
 (ii) \ $(f^\ell)^{-1}|_{L^k} = (f^k)^{-1}|_{L^k}$ \ ($\ell \geq k$). 
\end{itemize}
In fact, from $(2)_k$ and $(4)_k$ it follows that  
\[ \mbox{(a) $\phi^{k+1}|_{L_{k}} = id_{L_{k}}$ \ \ and \ \ 
(b) $\phi^{k+1}(B) = B$, 
$\ds {\rm dim}\,B \leq 2^{-k}$ $(B \in {\cal C}(L_k^c))$.} \]  
This implies that 
$\ds d(\phi^{k+1}, id_{M}) \leq 2^{-k}$ and 
$\ds d(f^{k}, f^{k+1}) \leq 2^{-k}$. 
On the other hand, \\ 
since (c) $(\phi^{k+1})^{-1}|_{L^{k}} = id_{L^{k}}$ and 
(d) $(\phi^{k+1})^{-1}(B) = B$ $(B \in {\cal C}(L_k^c))$, 
it follows that 
\[ \mbox{$(f^{k+1})^{-1}|_{L^{k}} = (f^{k})^{-1}|_{L^{k}}$ \ \ and \ \ 
$(f^{k+1})^{-1}(B) = (f^{k})^{-1}(B)$, \ \ whose diameter $\ds \leq 2^{-k}$ by $(4)_k$.} \] 
This implies 
$\ds d\big((f^{k})^{-1}, \big(f^{k+1})^{-1}\big) \leq 2^{-k}$. 
By $(2)_k$\,(ii), 
for $\ell \geq k+1$ we have   
$\phi^\ell|_{L^k} = id_{L^k}$ and 
\[ (f^\ell)^{-1}|_{L^k} = 
(\phi^{\ell} \cdots \phi^{k+1} f^k)^{-1}|_{L^k} 
= (f^k )^{-1} (\phi^{\ell} \cdots \phi^{k+1})^{-1}|_{L^k}
= (f^k)^{-1}|_{L^k}
. \] 

[2] There exists a canonical extension map $\eta : {\cal D}_\partial(M)_1 \to {\cal H}_{\partial \cup E_M}(\overline{M})$. By [1]\,(i)
the maps $\overline{f}^k = \eta f^k$ $(k \geq 1)$ satisfy the following conditions:  
\[ \mbox{$\ds d(\overline{f}^{k}, \overline{f}^{k+1}) \leq 2^{-k}$ \ \ and \ \ 
$\ds d\big(\big(\overline{f}^{k}\big)^{-1}, \big(\overline{f}^{k+1}\big)^{-1}\big) \leq 2^{-k}$ \ \ $(k \geq 1)$.} \] 
Since $\overline{M}$ is compact, 
the sequences $(\overline{f}^{k})_k$ and $\big((\overline{f}^{k})^{-1}\big)_k$ 
converge $d$-uniformly to some continuous maps 
$\overline{f}, \overline{f}' : P \to {\rm Map}_{\partial \cup E_M}(\overline{M})$. 
Since 
\[ \mbox{$\overline{f}' \overline{f} = \lim_{k\to\infty} (\overline{f}^{k})^{-1} \overline{f}^{k} \equiv id_{\overline{M}}$ \ \ and \ \ $\overline{f}\, \overline{f}' = \lim_{k\to\infty} \overline{f}^{k} (\overline{f}^{k})^{-1} \equiv id_{\overline{M}}$,} \]  
it follows that $\overline{f}_p \in {\cal H}_{\partial \cup E_M}(\overline{M})$ and 
$f_p = \overline{f}_p|_M \in {\cal H}_\partial(M)$ for each $p \in P$. 
Therefore the sequences $(f^k)_k$ and $\big((f^k)^{-1}\big)_k$ converge $d|_M$-uniformly to the maps $f = \overline{f}|_M, \ f^{-1} : P \to {\cal H}_\partial(M)$. 
In [1]\,(ii), tending $\ell \to \infty$, we have the condition (2) (i.e., 
$f^{-1}|_{L^k} = (f^{k})^{-1}|_{L^k}$ $(k \geq 1)$). 

[3] Since $(f^{k}_p)^{-1} \in {\cal D}_\partial(M)$ and 
$f_p^{-1}|_{L^k}$ is a $C^\infty$-embedding 
for any $k \geq 1$, we have $f_p^{-1} \in {\cal D}_\partial(M)$. 
To see that $f_p \in {\cal D}_\partial(M)_1$, 
we have to construct a path $F_p(t) \in {\cal D}_\partial(M)$ ($t \in [0, 1]$)
with $F_p(0) = id_M$ and $F_p(1) = f_p$. 
Since $\phi^k_p \in {\cal D}_{\partial \cup L^{k-1}}^c(M)_1^\ast$, 
\vspace{1mm} 
there exists a path 
$\Phi^k_p(t) \in {\cal D}_{\partial \cup L^{k-1}}^c(M)$ 
($k-1 \leq t \leq k$) such that 
$\Phi^k_p(k-1) = id_M$ and $\Phi^k_p(k) = \phi^k_p$. 
Replacing $[0,1]$ by $[0, \infty]$, we define $F_p(t) \in {\cal D}_\partial(M)$ ($t \in [0, \infty]$) by 
\[ F_p(t) = 
\begin{cases}
\Phi^k_p(t) f_p^{k-1} & (k \geq 1, \,k-1 \leq t \leq k) \\[1mm]
f_p & (t = \infty).  
\end{cases} 
 \] 
Obviously $F_p(t)$ is continuous on $[0, \infty)$. 
Since $F_p(t)^{-1}|_{L^k} = (f^{k}_p)^{-1}|_{L^k}$ ($k \leq t \leq \infty$), 
it follows that $F_p(t)^{-1} \in {\cal D}(M)$ is continuous at $t = \infty$.
Hence $F_p(t)$ is also contiuous at $t = \infty$. 

[4] Since the maps $(f^{k})^{-1} : P \to {\cal D}_\partial(M)_1$ ($k \geq 1$) are continuous, 
by the condition (2) and the definiton of the compact-open $C^\infty$-topology, 
it follows that the map $f^{-1} : P \to {\cal D}_\partial(M)_1$ is continuous. 
This implies the continuity of the map $f : P \to {\cal D}_\partial(M)_1$. 
This completes the proof.
\end{proof} 

Suppose $M$ is a noncompact, connected oriented $C^\infty$ $n$-manifold possibly with boundary 
and $L_i$ $(i \geq 0)$ is a sequence of compact $n$-submanifolds of $M$ 
such that 
\ $L_0 = \emptyset$, \ $\cup_i L_i = M$ \ and \ 
$L_{i-1} \subset {\rm Int}_M L_i$ \ $(i \geq 1)$. 
We set $S = \cup_i {\rm Fr}_M L_i$ and 
${\cal L} = \cup_i \,{\cal C}(cl(L_i - L_{i-1}))$. 

\begin{lemma}\label{lem_Moser} 
Suppose $\mu, \nu : P \to {\cal V}^+(M)$ are continuous maps 
such that $\mu(N) = \nu(N)$ $(N \in {\cal L})$. 
Then there exists a map $\chi : P \to {\cal D}_\partial(M, S)_1$ such that 
\[ \mbox{{\rm (1)} $\chi_\ast \mu = \nu$, \ \ 
{\rm (2)} $\chi(N) = N$ \ $(N \in {\cal L})$ \ \ and \ \ {\rm (3)} 
 if $p \in P$ and $\mu_p = \nu_p$, then $\chi_p = id_M$.} \]  
\end{lemma} 

\begin{proof} 
Choose a product collar $E^1 = \partial M \times [0,1]$ of $\partial M$ in $M$ and 
a product bicollar $E^2 = S \times [-1,1]$ of $S$ in $M$ such that 
$E^1$ intersects with $E^2$ orthogonally so that 
\[ \mbox{$E^1 \cap E^2 = \partial S \times [0,1] \times [-1,1]$ \quad and \quad 
$E^1[0,a] \cap E^2[-b, b] = \partial S \times [0,a] \times[-b, b]$ 
\quad $(a, b \in [0,1])$.} \] 

By Lemma~\ref{lem_collar} (2) there exist maps 
\[ \mbox{$\chi^1 : P \to {\cal D}_{\partial \cup (M \setminus E^1)}(M)_1$ 
\quad and \quad $\e^1 : P \to C^\infty(\partial M, (0,1))$} \] 
such that for each $p \in P$,  
\[ \mbox{(i) 
$(\chi^1_p)_\ast \mu_p = \nu_p$ on $E^1[0,\e^1_p]$,  \quad 
(ii) $\chi^1_p(E^1_x) = E^1_x$ $(x \in \partial M)$ \quad and \quad 
(iii) if $\mu_p = \nu_p$, then $\chi^1_p = id_M$.} \]
From (ii) it follows that $\chi^1_p \in {\cal D}_\partial(M; S)_1$ and 
$\chi^1_p(N) = N$ $(N \in {\cal L})$. 
Thus we obtain maps  
\[ \mbox{$\chi^1_p : P \to {\cal D}_\partial(M; S)_1$ \quad and \quad 
$\mu^1 = \chi^1_\ast \mu : P \to {\cal V}^+(M)$, \ 
$\mu^1_p = (\chi^1_p)_\ast \mu_p$.} \] 

By Lemma~\ref{lem_collar} (1) there exist maps 
\[ \mbox{$\chi^2 : P \to {\cal D}_{S \cup (M \setminus E^2)}(M)_1$ 
\quad and \quad 
$\e^2 : P \to C^\infty(S, (0,1))$} \]
such that for each $p \in P$, \quad (i)$'$ 
$(\chi^2_p)_\ast \mu^1_p = \nu_p$ on $E^2[-\e^2_p,\e^2_p]$, \quad and 
\begin{itemize} 
\item[(ii)$'$] for each $x \in S$ \ \ 
(a) $\chi^2_p((E^2)_x^\pm) =(E^2)_x^\pm$ \ and \ (b) 
if $\mu^1_p = \nu_p$ on $(E^2)_x^\pm$, 
then $\chi^2_p = id$ on $(E^2)_x^\pm$. 
\end{itemize}  
From (ii)$'$ it follows that $\chi^2_p \in {\cal D}_{\partial \cup S \cup (M \setminus E^2)}(M)_1$ and 
$\chi^2_p(N) = N$ $(N \in {\cal L})$. 
Thus we obtain maps  
\[ \mbox{$\chi^2_p : P \to {\cal D}_\partial(M; S)_1$ \quad and \quad 
$\mu^2 = \chi^2_\ast \mu^1 : P \to {\cal V}^+(M)$, \ 
$\mu^2_p = (\chi^2_p)_\ast \mu^1_p$.} \] 

For each $N \in {\cal L}$, let 
$\partial_- N = \partial N \cap \partial M$ and 
$\partial_+ N = {\rm Fr}_M N = \partial N \cap S$. 
Then  
$$E^N = (E^1 \cup E^2) \cap N = (E^2_{\partial_+ N} \cap N) \cup E^1_{\partial_- N}$$ 
is a collar of $\partial N$ in $N$. 
Consider the continuous maps $\mu^2|_N, \nu|_N : P \to {\cal V}^+(N)$ and 
$$\e^N : P \to (0,1/2), \quad 
\e^N_p = \frac{\,1}{\,2\,} \min \big\{
\min\,\e^1_p|_{\partial_- N}, \min \e^2_p|_{\partial_+ N} 
\big\}.$$ 
For each $p \in P$ it is seen that 
\[ \mbox{$\mu^2_p(N)  = \nu_p(N)$ \quad and \quad 
$\mu^2_p|_N = \nu_p|_N$ on $E^N[0, 2\e^N_p]$.} \]  
In fact, 
since $\chi^2_p(N) = \chi^1_p(N) = N$,  it follows that 
$\mu^2_p(N) = \mu^1_p(N) = \mu_p(N) = \nu_p(N)$. 
By (i)$'$ and (i) we have 
$\mu^2_p = \nu_p$ on $E^2_{\partial_+ N}[-2\e^N_p, 2\e^N_p]$  and 
$\mu^1_p = \nu_p$ on $E^1_{\partial_- N}[0, 2\e^N_p]$. 
The latter and (ii)$'$\,(b) imply that $\chi^2_p = id$ on $E^1_{\partial_- N}[0, 2\e^N_p]$ and so $\mu^2_p = \mu^1_p = \nu_p$ on $E^1_{\partial_- N}[0, 2\e^N_p]$. 

Then, Theorem~\ref{thm_Moser} yields  
a map $\chi^N : P \to {\cal D}_{\partial}(N)_1$ such that for each $p \in P$ 
\begin{itemize} 
\item[(i)$''$] $(\chi^N_p)_\ast (\mu^2_p|_N) = \nu_p|_N$, \ \  
(ii)$''$ $\chi^N_p = id$ on $E^N[0, \e^N_p]$ \ 
($\chi^N_p$ is isotopic to $id_N$ relative to $E^N[0, \e^N_p]$), \\[1mm] 
\hspace{-10mm} (iii)$''$ if $\mu^2_p|_N = \nu_p|_N$, then $\chi^N_p = id_N$.
\end{itemize} 
Define a map $\chi^3 : P \to {\cal D}_{\partial \cup S}(M)_1$ by 
$\chi^3_p|_N = \chi^N_p$ $(N \in {\cal L})$. 
Then we have $(\chi^3_p)_\ast \mu^2_p = \nu_p$ $(p \in P)$. 

Finally, the required map $\chi$ is defined by the composition 
\[ \mbox{
$\chi : P \to {\cal D}_\partial(M, S)_1$, \quad 
$\chi_p= \chi^3_p \,\chi^2_p \,\chi^1_p$ \ \ $(p \in P)$.} \] 
If $\mu_p = \nu_p$, then 
$\chi^1_p = id_M$ and $\mu^1_p = \nu_p$. In turn, 
this implies that $\chi^2_p = id_M$ and $\mu^2_p = \mu^1_p =\nu_p$ 
so that $\chi^3_p = id_M$ and hence $\chi_p = id_M$. 
\end{proof} 


\subsection{Main deformation theorem} \mbox{}

Suppose $M$ is a noncompact connected oriented $C^\infty$ $n$-manifold possibly with boundary, $d$ is a fixed metric on $\overline{M}$ 
and $\mu, \nu : P \to {\cal V}^+(M)$ are continuous maps. 
Let ${\cal F}(M)$ denote the class of connected $C^\infty$ $n$-submanifold $N$ of $M$ such that ${\rm Fr}_M N$ is compact 
and set ${\cal F}_c(M) = \{ N \in {\cal F}(M) \mid N : \text {compact}\}$. 
Let $C^0(P)$ denote the set of all continuous functions $\alpha : P \to {\Bbb R}$ 
and ${\cal D}$ denote the collection of all continuous maps 
$f : P \to {\cal D}_\partial(M)_1$. 

\begin{assumption}\label{asp_a} 
Suppose a subclass ${\cal F}$ of ${\cal F}(M)$ and 
a map $a : {\cal D}^2 \times {\cal F} \to C^0(P)$ 
satisfies the following conditions : 
\begin{itemize} 
\item[($\ast_0$)] 
${\cal F}_c(M) \subset {\cal F}$ and $\mu(A) = \nu(A) = \infty$ ($A \in {\cal F}(M) - {\cal F}$).  
\vskip 1mm 
\item[($\ast_1$)] $a(f, g; C) \in (- \mu(f^{-1}(C)), \nu(g^{-1}(C)))$ \ 
($f, g \in {\cal D}$, $C\in {\cal F}$).
\vskip 1mm 
\item[($\ast_2$)] $a(f_2, g; C) = a(f_1, g; C) + J^\mu(f_1^{-1}(C), f_2^{-1}(C))$ \\[1mm] 
$a(f, g_2; C) = a(f, g_1; C) - J^\nu(g_1^{-1}(C), g_2^{-1}(C))$ \hspace{3mm} \ \smash{\raisebox{4mm}{(\hspace{0.5mm}$f, f_1, f_2, \ g, g_1, g_2 \in {\cal D}$, \, $C \in {\cal F}$\hspace{0.5mm}).}}
\vskip 1mm 
\item[($\ast_3$)] If $(K, L) \in {\cal N}^{(2)}(M)$, $A \in {\cal C}(K^c)$ and 
${\cal C}(A \cap L^c) \subset {\cal F}$, then 
\[ \mbox{$A \in {\cal F}$ \ \ and  \ \ $a(f, g; A \cap L) + \sum_{\mbox{\tiny $B \!\in \!{\cal C}(A \cap L^c)$}} a(f, g; B) = a(f, g; A)$ \ \ ($f, g \in {\cal D}$). \hspace{10mm}} \] 
\item[($\ast_4$)] If $M \in {\cal F}$, then \ $a(id_M, id_M; M) = 0$. 
\end{itemize} 
\end{assumption} 

\begin{remark}\label{rem_a} From ($\ast_2$) it follows that for $p \in P$ 
\begin{itemize}
\item[(i)\,] $a_p(f_1, g_1; C) = a_p(f_2, g_2; C)$ \ if \ 
\begin{tabular}[t]{l}
$(f_1)_p^{-1}(C) = (f_2)_p^{-1}(C)$ \  and \ 
$(g_1)_p^{-1}(C) = (g_2)_p^{-1}(C)$ \\[2mm] 
\ \ (\,for example, $(f_1)_p = (f_2)_p$ \ and \ $(g_1)_p = (g_2)_p$\,), 
\end{tabular}
\vskip 1mm 
\item[(ii)] $a_p(\phi f, g; C) = a_p(f, g; C)$ \ if \ $\phi \in {\cal D}$ \ and \ $\phi_p = id$ on $C$. 
\end{itemize} 
\end{remark} 

See Section 4 for explicit examples of $(a, {\cal F})$. 

Recall the defining conditions $(1)_k \sim (4)_k$ for admissible sequences in Definition~\ref{def_adm_seq}. 

\begin{theorem}\label{thm_deform} {\bf [1]} There exists an admissible sequence $(K_k, L_k, f^k, g^k)$ $(k=1,2,\cdots)$ $(L^0 = \emptyset$, $f^0 = g^0 \equiv id_M)$ which satisfies the following additional conditions:  
\begin{itemize} 
\item[$(5)_k$]\hspace{-2mm} 
\begin{tabular}[t]{cllll}
{\rm (i)} & $a(f^k, g^{k-1}; C) = 0$ & $\big(\, C \in {\cal C}(cl(K_k - L_{k-1})) \cup ({\cal C}(K_k^c) \cap {\cal F}))\, \big)$ \\[1mm] 
{\rm (ii)} & $a(f^k, g^k; C) = 0$ & $\big(\, C \in {\cal C}(cl(L_k - K_k)) \cup  ({\cal C}(L_k^c) \cap {\cal F})) \, \big)$
\end{tabular} 
\vskip 1mm
\item[$(6)_k$]  If $p \in P$ 
and $a_p(id_M, id_M; C) = 0$ $(C \in {\cal F})$, then $f^k_p = g^k_p = id_M$. 
\end{itemize} 
{\bf [2]} The sequences $(f^k)_k$, $(g^k)_k$ 
converge $d|_M$-uniformly to some maps $f, g : P \to {\cal D}_\partial(M)_1$ respectively. The maps $f$ and $g$ satisfy the following conditions: 
\begin{itemize}
\item[(1)] $f^{-1}|_{L_k} = (f^k)^{-1}|_{L_k}$, \ $g^{-1}|_{K_k} = (g^{k-1})^{-1}|_{K_k}$ \ $(k \geq 1)$. 
\vskip 1mm
\item[(2)] $a(f, g; C) = 0$ for any $C \in {\cal F}$ with \\ 
\hspace{15mm} $C \, \in \, {\cal C}(K_k^c) \cup {\cal C}(L_k^c) \cup 
{\cal C}(cl(K_k - L_{k-1})) \cup {\cal C}(cl(L_k - K_k))$ \ for some \ $k \geq 1$.  
\vskip 1mm
\item[$(3)$]  If $p \in P$ and $a_p(id_M, id_M; C) = 0$ $(C \in {\cal F})$, then $f_p = g_p = id_M$. 
\end{itemize} 
\end{theorem} 

\begin{proof} 
{\bf [1]} Suppose $k \geq 1$ and assume that we have constructed $(K_i, L_i, f^i, g^i)$ ($i = 1, \cdots, k-1$) which satisfy the conditions $(1)_i \sim (6)_i$ ($i = 1, \cdots, k-1$).  
The tuple $(K_k, L_k, f^k, g^k)$ is constructed as follows: 

$K_k$ :  Since $\big\{ (g^{k-1}_p)^{-1} \big\}_p$ is equi-continuous with respect to $d|_M$, 
we can find $K_k \in {\cal N}(M)$ which satisfies $(1)_k$ and $(4)_k$\,(i). 

$f^k$ : 
For each $B \in {\cal C}(L_{k-1}^c)$ we construct a map 
$\phi^B : P \to {\cal D}^c_{\partial \cup B^c}(M)_1^\ast$ which satisfies the following conditions: 
\begin{itemize} 
\item[(3)$_k'$] 
\begin{itemize} 
\item[(i)\,] $\phi^B$ has locally common compact support in ${\rm Int}\,M$. 
\item[(ii)] $\{ (\phi^B_p)^{-1} \}_p$ is equicontinuous with respect to $d|_M$. 
\end{itemize} 
\vskip 1mm
\item[(5)$_k'$] $a(\phi^B f^{k-1}, g^{k-1}; B \cap K_k) = 0$ \ and \ 
$a(\phi^B f^{k-1}, g^{k-1}; A) = 0$ $(A \in {\cal C}(B \cap K_k^c) \cap {\cal F})$. \vskip 1mm
\item[(6)$_k'$]  If $p \in P$ and $a_p(id_M, id_M; C) = 0$ $(C \in {\cal F})$, then $\phi^B_p = id_M$. 
\end{itemize} 
Consider the decomposition of $B$ and the associated functions in $C^0(P)$ defined by 
\begin{itemize}
\item[(i)\,] $B = N \cup \big( \cup_{i=1}^\ell A_i\big) \cup \big( \cup_{j=1}^m A_j' \big)$ \hspace{10mm} $(\ell, m \geq 0, \ \ell + m \geq 1)$,  \\[1mm] 
where 
$N = A_0 = B \cap K_k$, \ ${\cal C}(B \cap K_k^c) \cap {\cal F} = \{ A_1, \cdots, A_\ell \}$, \ ${\cal C}(B \cap K_k^c) - {\cal F} = \{ A_1', \cdots, A_m'\}$ 
\vskip 2mm
\item[(ii)] 
$a(i) = a(f^{k-1}, g^{k-1}; A_i)$ \ ($i = 0, 1, \cdots, \ell$), \ \ 
$\ds b(j) = - \frac{1}{m} \mbox{$\sum_{i=0}^\ell a(i)$}$ \ ($j = 1, \cdots, m$). 
\end{itemize}
\vskip 1mm 
In order to apply Lemma~\ref{lem_deform} to these data and $(f^{k-1}_{\ \ \  \ast} \mu)|_B$, 
we have to check the next conditions: 
\begin{itemize}
\item[(a)] $a(0) + \sum_{i=1}^\ell a(i) + \sum_{j=1}^m b(j) = 0$, 
\vskip 2mm 
\item[(b)] $a(i) > -\mu((f^{k-1})^{-1}(A_i))$ \ ($i = 0, 1, \cdots, \ell$), \ \ 
$b(j) > -\mu((f^{k-1})^{-1}(A_j'))$ \ ($j = 1, \cdots, m$). 
\end{itemize}

For (a), if $m \geq 1$, then this follows from the definition of $b(j)$'s. 
If $m = 0$, then from $(\ast_3)$ and  $(5)_{k-1}$\,(ii) ($(\ast_4)$ for $k = 1$) it follows that \ 
$B \in{\cal F}$ \ and \hspace{2mm} 
$\sum_{i=0}^\ell a(i) = a(f^{k-1}, g^{k-1}; B) = 0$.

For (b), the assertion for $a(i)$'s follows from $(\ast_1)$. 
Since $A_j' \not\in {\cal F}$, from $(\ast_0)$ it follows that 
$\mu(A_j') = \infty$. Since each $f^{k-1}_p$ has a compact support, we also have 
$\mu((f^{k-1})^{-1}(A_j')) = \infty$, which obviously implies the  
assertion for $b(j)$'s. 

\vskip 1mm 
Now Lemma~\ref{lem_deform} can be applied to yield a map 
$\phi^B : P \to {\cal D}^c_{\partial \cup B^c}(M)_1^\ast$ which satisfies $(3)_k'$ and the next conditions: 
\vskip 1mm 
\begin{itemize} 
\item[(5)$_k''$]  
$J^{f^{k-1}_{\ \ \, \ast} \mu}((\phi^B)^{-1}(A_i), A_i) = a(i)$ ($i = 0, 1, \cdots, \ell$), \ 
$J^{f^{k-1}_{\ \ \, \ast} \mu}((\phi^B)^{-1}(A_j'), A_j') = b(j)$ ($j = 1, \cdots, m$). 
\vskip 1mm 
\item[(6)$_k''$]  If $p \in P$ and $a_p(i) = 0$ ($i = 0,1, \cdots, \ell$) (so that $b_p(j) = 0$ ($j = 1, \cdots, m$)), then $\phi^B_p = id_M$. 
\end{itemize} 
\vskip 1mm 
We note that the conditions (5)$_k''$ and (6)$_k''$ imply 
(5)$_k'$ and (6)$_k'$ respectively. 
In fact, by $(\ast_2)$ and Lemma~\ref{lem-difference of volume} 
the condition $(5)_k''$ implies that 
\begin{align*}
a(\phi^Bf^{k-1}, g^{k-1}; A_i) &= 
a(f^{k-1}, g^{k-1}; A_i) + J^\mu((f^{k-1})^{-1}(A_i), (\phi^Bf^{k-1})^{-1}(A_i)) \\[1mm] 
&= a(i) - J^{f^{k-1}_{\ \ \, \ast} \mu}((\phi^B)^{-1}(A_i), A_i) = 0 \hspace{10mm}  
(i = 0, 1, \cdots, \ell). 
\end{align*}
If $p \in P$ and $a_p(id_M, id_M; C) = 0$ $(C \in {\cal F})$, then 
from $(6)_{k-1}$ and Remark 3.3\,(i) it follows that 
\[ \mbox{$f^{k-1}_p = g^{k-1}_p = id_M$ \ \ and \ \ 
$a_p(i) = a_p(f^{k-1}, g^{k-1}; A_i) = a_p(id_M, id_M; A_i) = 0 \hspace{5mm} 
(i = 0, 1, \cdots, \ell)$,} \]
and from (6)$_k''$ follows $\phi^B_p = id_M$. 

Finally, we define a map $\phi^k : P \to {\cal D}_{\partial \cup L_{k-1}}^c(M)_1^\ast$ by 
\[ \mbox{$\phi^k_p = \phi^{B_1}_p \cdots \phi^{B_r}_p$ \ ($p \in P$), \ where \ 
${\cal C}(L_{k-1}^c) = \{ B_1, \cdots, B_r \}$ \ is any ordering.} \] 
Then the required map $f^k : P \to {\cal D}_\partial^c(M)_1^\ast$ is defined by $f^k = \phi^k f^{k-1}$. 
The conditions $(3)_k$ and $(6)_k$ follow from 
$(3)_{k-1}$,  $(3)_k'$ and $(6)_{k-1}$,  $(6)_k'$ respectively, while  
the condition $(5)_k$ follows from $(5)_k'$ and Remark~\ref{rem_a} (ii). 
\vskip 1mm 
$L_k$ : Since $\big\{ (f^k_p)^{-1} \big\}_p$ is equi-continuous, 
we can find $L_k \in {\cal N}(M)$ which satisfies $(1)_k$ and $(4)_k$\,(ii). 

$g^k$ : The map $\psi^k$ and $g^k$ are obtained by the same arguments, using Lemma~\ref{lem_deform} and Remark~\ref{rem_deform}\,(B). 
\vskip 1mm 
{\bf [2]} The existence of the limit maps $f$ and $g$ and the condition (1) follow from Lemma~\ref{lem_lim}. 
\begin{itemize} 
\item[(2)] For each $k \geq 1$ 
\begin{itemize} 
\item[(i)\,] since $f^{-1} = (f^k)^{-1}$ and $g^{-1} = (g^{k-1})^{-1}$ on $K_k$, 
for any $C \in ({\cal C}(K_k^c) \cap {\cal F}) \cup {\cal C}(cl(K_k - L_{k-1}))$ it follows that \ \ $f^{-1}(C) = (f^k)^{-1}(C)$, $g^{-1}(C) =(g^{k-1})^{-1}(C)$ \ \ and hence \\
\hspace{10mm} $a(f, g; C) = a(f^{k}, g^{k-1}; C) = 0$ \ \ 
by Remark~\ref{rem_a}\,(i) and  $(5)_k$\,(i), 
\vskip 1mm
\item[(ii)]  since $f^{-1} = (f^k)^{-1}$ and $g^{-1} = (g^k)^{-1}$ on $L_k$, 
for any $C \in ({\cal C}(L_k^c) \cap {\cal F}) \cup {\cal C}(cl(L_k - K_k))$ it follows that \ \ $f^{-1}(C) = (f^k)^{-1}(C)$, $g^{-1}(C) =(g^k)^{-1}(C)$ \ \ and hence \\
\hspace{10mm} $a(f, g; C) = a(f^{k}, g^k; C) = 0$ \ \ 
by Remark~\ref{rem_a}\,(i) and  $(5)_k$\,(ii). 
\end{itemize}
\vskip 1mm
\item[(3)]  The assertion follows from $(6)_k$. 
\end{itemize} 
\vskip -7mm 
\end{proof} 

\section{Proof of Theorems 1.1 and 1.2} 

In this section we apply Theorem~\ref{thm_deform} to 
deduce Theorems 1.1 and 1.2. 
Suppose $M$ is a noncompact connected oriented $C^\infty$ $n$-manifold 
possibly with boundary 
and $d$ is a fixed metric on $\overline{M}$. 

\subsection{Proof of Theorem 1.1} \mbox{} 

Suppose $F$ is an open subset of $E_M$ 
and 
$\mu, \nu : P \to {\cal V}^+(M; F)_{ew}$ are continuous maps with $\mu_p(M) = \nu_p(M)$ $(p \in P)$. 

\begin{defi} We define $({\cal F}, a)$ as follows : 
\begin{itemize}
\item[(i)\,] ${\cal F} = \big\{ C \in {\cal F}(M) \mid E_C \subset F \big\}$,  
\item[(ii)] $a : {\cal D}^2 \times {\cal F} \to C^0(P)$ : \ \ 
$a(f, g; C) = \nu(g^{-1}(C)) - \mu(f^{-1}(C)) = (g_\ast\nu)(C) - (f_\ast\mu)(C)$. 
\end{itemize}
\end{defi} 

\begin{remark}\label{rem_f} For $C \in {\cal F}(M)$ and $f \in {\cal D}$ we have  
\begin{itemize}
\item[(a)] $C \in {\cal F}$ iff $\mu(C) < \infty$ (or $\nu(C) < \infty$), 
since $E_M^{\mu_p} = E_M^{\nu_p} = F$ $(p \in P)$. 
\item[(b)] $E_{f^{-1}(C)} = E_C$ and so $C \in {\cal F}$ iff $f^{-1}(C) \in {\cal F}$, 
since $f_p \in {\cal D}_{E_M}(M)$ ($p \in P$).   
\end{itemize}
\end{remark}

\begin{lemma} The map $a$ is well-defined and $({\cal F}, a)$ satisfies the conditions $(\ast_0) \sim (\ast_4)$ in Assumption~\mbox{\rm \ref{asp_a}}. 
\end{lemma} 

\begin{proof} 
Remark~\ref{rem_f} implies that $\mu(f^{-1}(C)), \nu(g^{-1}(C)) < \infty$. 
Since $\mu(f^{-1}(C)) = J^\mu(f^{-1}(C), C) + \mu(C)$, 
by Lemmas~\ref{lem_conti} and \ref{lem_conti_ew} the map 
$\mu(f^{-1}(C)) : P \to {\Bbb R}$ is continuous. 
Similarly, the map 
$\nu(g^{-1}(C)) : P \to {\Bbb R}$ is continuous. 
Thus, the map $a(f,g;C) : P \to {\Bbb R}$ is continuous and 
the map $a$ is well-defined. 
The conditions $(\ast_0) \sim (\ast_4)$ are verified as follows. 

($\ast_0$) Obviously ${\cal F}_c(M) \subset {\cal F}$ and 
by Remark~\ref{rem_f}\,(a) 
we have $\mu(A) = \nu(A) = \infty$ ($A \in {\cal F}(M) - {\cal F}$).  
\vskip 1mm 
\begin{itemize} 
\item[($\ast_1$)] $a(f, g; C) = \nu(g^{-1}(C)) - \mu(f^{-1}(C)) 
\in (- \mu(f^{-1}(C)), \nu(g^{-1}(C)))$. 
\vskip 1mm 
\item[($\ast_2$)] Since $\mu(f_i^{-1}(C)), \nu(g_i^{-1}(C)) < \infty$ ($i=1,2$), we have 
\begin{align*}
a(f_1, g; C) + J^\mu(f_1^{-1}(C), f_2^{-1}(C)) 
&= \nu(g^{-1}(C)) - \mu(f_1^{-1}(C)) + \mu(f_1^{-1}(C)) - \mu(f_2^{-1}(C)) \\[1mm] 
&= \nu(g^{-1}(C)) - \mu(f_2^{-1}(C)) 
= a(f_2, g; C), \\[1mm]
a(f, g_1; C) - J^\nu(g_1^{-1}(C), g_2^{-1}(C)) 
&= \nu(g_1^{-1}(C)) - \mu(f^{-1}(C)) - \big(\nu(g_1^{-1}(C)) - \nu(g_2^{-1}(C))\big) \\[1mm] 
&= \nu(g_2^{-1}(C)) - \mu(f^{-1}(C)) = a(f, g_2; C).
\end{align*}
\vskip 1mm 
\item[($\ast_3$)] If $(K, L) \in {\cal N}^{(2)}(M)$, $A \in {\cal C}(K^c)$ and 
${\cal C}(A \cap L^c) \subset {\cal F}$, then 
\begin{itemize} 
\item[(i)\,] $\mu(A) = \mu(A \cap L) + \sum_{\mbox{\tiny $B \!\in \!{\cal C}(A \cap L^c)$}}  \mu(B) < \infty$, so that $A \in {\cal F}$, 
\vskip 1mm
\item[(ii)] $a(f, g; A \cap L) + \sum_{\mbox{\tiny $B \!\in \!{\cal C}(A \cap L^c)$}} a(f, g; B) \\[1mm]
\hspace{10mm} = \nu(g^{-1}(A \cap L)) - \mu(f^{-1}(A \cap L))
+ \sum_{\mbox{\tiny $B \!\in \!{\cal C}(A \cap L^c)$}}
\big(\nu(g^{-1}(B)) - \mu(f^{-1}(B))\big) \\[1mm] 
\hspace{10mm} =  \nu(g^{-1}(A)) - \mu(f^{-1}(A)) = a(f, g; A).$ 
\end{itemize} 
\vskip 1mm 
\item[($\ast_4$)] If $M \in {\cal F}$, then  
$a(id_M, id_M; M) = \nu(M) - \mu(M) = 0$ by the assumption on $\mu$ and $\nu$. 
\end{itemize} 
\vskip -7mm
\end{proof} 

\begin{proof}[\bf Proof of Theorem 1.1]  
By Theorem~\ref{thm_deform} 
we obtain an admissible sequence $(K_k, L_k, f^k, g^k)_{k \geq 1}$ and the limit maps  
$f, g : P \to {\cal D}_\partial(M)_1$ which satisfy the following conditions: 
\begin{itemize}
\item[(1)] $a(f, g; C) = 0$ (i.e., $(f_\ast\mu)(C) = (g_\ast\nu)(C)$) \\
\hspace{50mm} for any \ $C \in {\cal C}(cl(K_k - L_{k-1})) \cup {\cal C}(cl(L_k - K_k))$ \ ($k \geq 1$).
\vskip 1mm
\item[$(2)$]  If $p \in P$ and $a_p(id_M, id_M; C) = 0$ ($C \in {\cal F}$) (i.e., $\mu_p = \nu_p$), then $f_p = g_p = id_M$. 
\end{itemize} 
By the condition (1) and Lemma~\ref{lem_Moser} we obtain a map 
$\chi : P \to {\cal D}_\partial(M)_1$ such that \\
\hspace{10mm} (i) $\chi_\ast (f_\ast\mu) = g_\ast\nu$ \ \ and \ \ (ii) 
if $p \in P$ and $(f_p)_\ast\mu_p = (g_p)_\ast\nu_p$, then $\chi_p = id_M$. \\
Finally, the required map $h :  P \to {\cal D}_\partial(M)_1$ is defined by \ \ $h_p = g_p^{-1}\chi_p f_p$ \ ($p \in P$). 
\end{proof} 

Theorem 1.1 is equivalent to the following statement. 

\begin{theorem-v-1}\label{thm_v-1}
Suppose $P$ is any topological space and 
$\widetilde{\mu} : P \times [0,1] \to {\cal V}^+(M; F)_{ew}$ is a homotopy with $\widetilde{\mu}_p(t)(M) = \widetilde{\mu}_p(0)(M)$ $((p,t) \in P \times[0,1])$. 
Then there exists a homotopy $\widetilde{h} : P \times [0,1] \to {\cal D}_\partial(M)_1$ such that for each $(p,t) \in P \times[0,1]$ 
\[ \mbox{{\rm (1)} $\widetilde{h}_p(t)_\ast \widetilde{\mu}_p(0) = \widetilde{\mu}_p(t)$, \ {\rm (2)} $\widetilde{h}_p(0) = id_M$  \ and \ {\rm (3)} if $\widetilde{\mu}_p(t) = \widetilde{\mu}_p(0)$, then $\widetilde{h}_p(t) = id_M$.} \]  
\end{theorem-v-1}

\begin{proof} 
The homotopy $\widetilde{h}$ in Theorem 1.1$'$ 
is obtained by applying Thereom 1.1 to the maps $\widetilde{\mu}$ 
and $\widetilde{\mu}' : P \times [0,1] \to {\cal V}^+(M; F)_{ew}$, $\widetilde{\mu}'(t) = \widetilde{\mu}(0)$ $(t \in [0,1])$. 

The maps $\mu, \nu : P \to {\cal V}^+(M; F)_{ew}$ in Theorem 1.1 
are connected by the affine homotopy \break 
$\widetilde{\mu} : P \times [0,1] \to {\cal V}^+(M; F)_{ew}$, $\widetilde{\mu}(t) = (1-t)\mu + t \nu$. 
Theorem 1.1$'$ yields a homotopy $\widetilde{h}$ which traces $\widetilde{\mu}$,  and the map $h = \widetilde{h}(1)$ satisfies the required conditions in Theorem 1.1. 
\end{proof}

\begin{proof}[\bf Proof of Corollary 1.1\,(1)] 
We apply Theorem 1.1 to the data: \\  
$P = {\cal V}^+(M; m, F)$, \ 
the constant map $P \to {\cal V}^+(M; F)$, $\mu \mapsto \omega$ \ and \ 
the inclusion map $P \subset {\cal V}^+(M; F)$. \\[1mm]  
This yields a map $\sigma : P \to {\cal D}(M)_1$ such that 
$\sigma(\mu)_\ast \omega = \mu$ ($\mu \in P$) and $\sigma(\omega) = id_M$. 
This means that the map $\sigma$ is the required section of $\pi_\omega$. 
\end{proof} 

Corollary 1.1\ (2) is a special case of the next corollary. 

\begin{cor}\label{cor_sdr} Suppose $G$ is any subgroup of ${\cal D}^+(M; F)$ 
with ${\cal D}_\partial(M)_1 \subset G$. Then 
\begin{itemize}
\item[(1)] $(G, G \cap {\cal D}(M; \omega)) \cong 
({\cal V}^+(M; m, F)_{ew}, \{ \omega \}) \times (G \cap {\cal D}(M; \omega))$.
\item[(2)] $G \cap {\cal D}(M; \omega)$ is a strong deformation retract of $G$. 
\end{itemize} 
\end{cor}

\begin{proof} 
(1) The required homeomorphism  
$\phi : G \to {\cal V}^+(M; m, F)_{ew} \times (G \cap {\cal D}(M; \omega))$ is defined by 
$\phi(h) = (\pi_\omega(h), (\sigma(\pi_\omega(h)))^{-1}h)$. 
The inverse is given by $\phi^{-1}(\mu, h) = \sigma(\mu)h$. 

(2) Since the singleton $\{ \omega \}$ is a strong deformation retract of ${\cal V}^+(M; m, F)_{ew}$, 
the conclusion follows from (1). 
\end{proof}

\begin{remark}\label{rem_sdr} 
${\cal D}_\partial(M)_1 \cap {\cal D}(M; \omega) 
= {\cal D}_\partial(M; \omega)_1$. 
\end{remark} 

\begin{proof}
By Corollary~\ref{cor_sdr}\,(2) ${\cal D}_\partial(M)_1 \cap {\cal D}(M; \omega)$ is a strong deformation retract of ${\cal D}_\partial(M)_1$. Since the latter is path-connected, so is the former. 
\end{proof} 

\subsection{Proof of Theorem 1.2} \mbox{} 

Suppose $\mu : P \to {\cal V}^+(M)$ and $a : P \to {\cal S}(M)$ are continuous maps such that $a_p \in {\cal S}(M, \mu_p)$ $(p \in P)$. 
(Let $\nu = \mu$.) 

\begin{defi} We define $({\cal F}, \widetilde{a})$ as follows : 
\begin{itemize}
\item[(i)\,] ${\cal F} = {\cal F}(M)$,  
\item[(ii)] $\widetilde{a} : {\cal D}^2 \times {\cal F} \to C^0(P)$ : \ \ 
$\widetilde{a}(f, g; C) = a(E_C) - J^\mu(f^{-1}(C), g^{-1}(C))$. 
\end{itemize}
\end{defi} 

\begin{lemma} The map $\widetilde{a}$ is well-defined and $({\cal F}, \widetilde{a})$ satisfies the conditions $(\ast_0) \sim (\ast_4)$ in Assumption~\mbox{\rm \ref{asp_a}}. 
\end{lemma} 

The map $a$ is well-defined and $({\cal F}, a)$ satisfies the conditions $(\ast_0) \sim (\ast_4)$ in Assumption~\mbox{\rm \ref{asp_a}}. 

\begin{proof} 
Lemma~\ref{lem_conti} implies that the function 
$J^\mu(f^{-1}(C), g^{-1}(C)) : P \to {\Bbb R}$ is well-defined and continuous, 
while the continuity of the function $a(E_C) : P \to {\Bbb R}$ follows from 
the definition of the weak topology of ${\cal S}(M)$. 
The conditions $(\ast_0) \sim (\ast_4)$ are verified as follows. 
\begin{itemize}
\item[($\ast_0$)] The assertion is trivial. 
\vskip 1mm 
\item[($\ast_1$)] Let $(f, g, C) \in {\cal D}^2 \times {\cal F}$ and $p \in P$. 
Since $f_p, g_p \in {\cal D}_{E_M}(M)$, we have  
$E_C = E_{f_p^{-1}(C)} = E_{g_p^{-1}(C)}$. 
If $E_C \subset E_M^{\mu_p}$, then (a) $a_p(E_C) = 0$ since $a_p \in {\cal S}(M; \mu_p)$ and (b) $\mu_p(f_p^{-1}(C)), \mu_p(g_p^{-1}(C)) < \infty$, so that 
\[ \widetilde{a}_p(f, g; C) = -J^{\mu_p}(f_p^{-1}(C), g_p^{-1}(C)) 
= \mu_p(g_p^{-1}(C)) - \mu_p(f_p^{-1}(C)). \] 
If $E_C \not\subset E_M^{\mu_p}$, then 
$\mu_p(f_p^{-1}(C)) = \mu_p(g_p^{-1}(C)) = \infty$. \\ 
In both cases, we have $\widetilde{a}_p(f, g; C) \in (- \mu_p(f_p^{-1}(C)), \mu_p(g_p^{-1}(C)))$. 
\vskip 1mm 

\item[($\ast_2$)]  By Lemma~\ref{lem-difference of volume} it follows that 
\begin{align*}
\mbox{} \hspace{10mm} 
\widetilde{a}(f_1, g; C) + J^\mu(f_1^{-1}(C), f_2^{-1}(C)) 
&= a(E_C) - J^\mu(f_1^{-1}(C), g^{-1}(C)) - J^\mu(f_2^{-1}(C), f_1^{-1}(C)) \\[1mm]  
&= a(E_C) - J^\mu(f_2^{-1}(C), g^{-1}(C)) = \widetilde{a}(f_2, g; C), \\[1mm]
\widetilde{a}(f, g_1; C) - J^\mu(g_1^{-1}(C), g_2^{-1}(C)) 
&= a(E_C) - J^\mu(f^{-1}(C), g_1^{-1}(C)) - J^\mu(g_1^{-1}(C), g_2^{-1}(C)) \\[1mm] 
&= a(E_C) - J^\mu(f^{-1}(C), g_2^{-1}(C)) = \widetilde{a}(f, g_2; C). 
\end{align*}
\vskip 1mm 

\item[($\ast_3$)] Let $(K, L) \in {\cal N}^{(2)}(M)$, $A \in {\cal C}(K^c)$ and $f, g \in {\cal D}$. Set $N = A \cap L$ and ${\cal C}(A \cap L^c) = \{ B_1, \cdots, B_m\}$. 
Note that (a) $a(E_N) = 0$ since $N$ is compact and $E_N = \emptyset$ and 
(b) $a(E_A) = \sum_i a(E_{B_i})$ since 
$E_A$ is the disjoint union of $E_{B_i}$ ($i=1, \cdots, m$). 
Thus by Lemma~\ref{lem-difference of volume}\,(3) we have 
\begin{align*}
\mbox{} \hspace{15mm} \widetilde{a}(f, g; N) + \mbox{$\sum_i$ } \widetilde{a}(f, g; B_i) 
&= - J^\mu(f^{-1}(N), g^{-1}(N)) 
+ \mbox{$\sum_i$}\,\big(a(E_{B_i}) - J^\mu(f^{-1}(B_i), g^{-1}(B_i))\big) \\[1mm]  
&= a(E_A) - J^\mu(f^{-1}(A), g^{-1}(A)) = \widetilde{a}(f, g; A).
\end{align*}

\item[($\ast_4$)] $\widetilde{a}(id_M, id_M; M) = a(E_M) - J^\mu(M, M) = 0$. 
\end{itemize}
\vskip -7mm 
\end{proof} 

\begin{proof}[\bf Proof of Theorem 1.2]  
Theorem~\ref{thm_deform} yields an admissible sequence $(K_k, L_k, f^k, g^k)_{k \geq 1}$ and the limit maps $f, g : P \to {\cal D}_\partial(M)_1$ 
which satisfy the following conditions: 
\begin{itemize}
\item[(1)] $\widetilde{a}(f, g; C) = 0$ 
\hspace{2mm} for any \ 
$C \, \in \, {\cal C}(cl(K_k - L_{k-1})) \cup {\cal C}(cl(L_k - K_k)) \cup {\cal C}(K_k^c) \cup {\cal C}(L_k^c)$ \ ($k \geq 1$). 
\vskip 1mm
\item[$(2)$]  If $p \in P$ and $\widetilde{a}_p(id_M, id_M; C) = 0$ ($C \in {\cal F}$), then $f_p = g_p = id_M$. 
\end{itemize} 

For any $C \, \in \, {\cal C}(cl(K_k - L_{k-1})) \cup {\cal C}(cl(L_k - K_k))$ ($k \geq 1$), since $C$ is compact, it follows that 
$E_C = \emptyset$ and $\mu(f^{-1}(C)), \mu(g^{-1}(C)) < \infty$, so that 
$$0 = \widetilde{a}(f, g; C) = -J^{\mu}(f^{-1}(C), g^{-1}(C)) 
= (g_\ast \mu)(C) - (f_\ast \mu)(C).$$ 
Thus, by Lemma~\ref{lem_Moser} we obtain a map 
$\chi : P \to {\cal D}_\partial(M)_1$ such that 
\begin{itemize}
\item[(3)] \ (i) $\chi_\ast (f_\ast\mu) = g_\ast\mu$, \ \ 
(ii) if $p \in P$ and $(f_p)_\ast\mu_p = (g_p)_\ast\mu_p$, then $\chi_p = id_M$ \ \ \ and \\ 
(iii) $\chi(C) = C$
\ ($C \in {\cal C}(cl(K_k - L_{k-1})) \cup {\cal C}(cl(L_k - K_k)), k \geq 1$). 
\end{itemize} 

We show that the map \ $h :  P \to {\cal D}_\partial(M)_1$, \ $h_p = g_p^{-1}\chi_p f_p$ \ ($p \in P$) \ satisfies the required conditions. 
\begin{itemize}
\item[(a)] By Remark~\ref{rem_sdr} we have 
$h_p \in {\cal D}_\partial(M)_1 \cap {\cal D}_\partial(M; \mu_p) 
= {\cal D}_\partial(M; \mu_p)_1.$ 
\item[(b)] For any $A \in {\cal C}(K_k^c)$ ($k \geq 1$), 
since $f_p \in {\cal D}_{E_M}(M)$ ($p \in P$), we have $E_A = E_{f^{-1}(A)}$, 
and by the condition (1) it follows that  
$$c^\mu_h(E_A) = c^\mu_h(E_{f^{-1}(A)}) 
= J^\mu(f^{-1}(A), h f^{-1}(A)) = J^\mu(f^{-1}(A), g^{-1}(A)) = a(E_A).$$ 
In general, for any $F \in {\cal Q}(E_M)$, 
there exists $k \geq 1$ and $A_1, \cdots, A_\ell \in {\cal C}(K_k^c)$ such that 
$F = \cup_i E_{A_i}$ (a disjoint union). 
Then, we have \ \ $c^\mu_h(F) = \sum_i c^\mu_h(E_{A_i}) = \sum_i a(E_{A_i}) = a(F)$. \\  
This implies that \ $c^\mu_h = a$.
\item[(c)] If $a_p = 0$, then $\widetilde{a}_p(id_M, id_M; C) = 0$ ($C \in {\cal F}$). Thus $f_p = g_p = id_M$ by the condition (2) and so $\chi_p = id_M$ by the condition (3)\,(ii). 
This implies that $h_p = id_M$. 
\end{itemize} 
This completes the proof. 
\end{proof} 

Theorem~\ref{thm_a} has the following slight generalization. 

\begin{theorem-a-1}\label{thm_a-1} 
Suppose $P$ is any topological space and 
$\mu : P \to {\cal V}^+(M)$ and $a : P \to {\cal S}(M)$ are continuous maps such that $a_p \in {\cal S}(M; \mu_p)$ $(p \in P)$. 
Then there exists a homotopy $h : P \times [0,1] \to {\cal D}_\partial(M)_1$ 
such that for each $(p, t) \in P \times [0,1]$ \\[1mm]  
\hspace{2mm} {\rm (1)} $h_p(t) \in {\cal D}_\partial(M; \mu_p)_1$, \ 
{\rm (2)} $h_p(0) = id_M$, \ 
{\rm (3)} $c_{h_p(t)}^{\mu_p} = t a_p$ \ and \ 
{\rm (4)} if $a_p = 0$, then $h_p(t) = id_M$.  
\end{theorem-a-1}

\begin{proof}
The homotopy $h$ is obtained by applying Theorem~\ref{thm_a} to the maps 
\[ \mbox{$\overline{\mu} : P \times [0,1] \to {\cal V}^+(M)$, \ $\overline{\mu}_{(p,t)} = \mu_p$ \ \ and \ \ 
$\overline{a} : P \times [0,1] \to {\cal S}(M)$, \ $\overline{a}_{(p,t)} = ta_p$.} \] 
\vskip -7mm  
\end{proof} 

\begin{proof}[\bf Proof of Corollary 1.2\,(1)] 
The required section is obtained by applying Theorem 1.2 to the data: 
$P = {\cal S}(M; \omega)$, the constant map 
$\mu \equiv \omega : P \to {\cal V}^+(M)$ and 
$id : P \to {\cal S}(M; \omega)$. 
\end{proof} 

Suppose $G$ is any subgroup of ${\cal D}_{E_M}(M; \omega)$ with 
${\cal D}_\partial(M; \omega)_1 \subset G$. 
Consider the restriction $c^\omega|_G : G \to S(M; \omega)$. 
Corollary 1.2\ (2) is a special case of the next corollary. 

\begin{cor} {\rm (1)} $(G, {\rm ker}\,c^\omega|_G) \cong ({\cal S}(M; \omega), 0) \times {\rm ker}\,c^\omega|_G$. 
\begin{itemize}
\item[(2)] ${\rm ker}\,c^\omega|_G$ is a strong deformation retract of $G$.
\end{itemize}
\end{cor} 

\begin{proof} 
(1) The required homeomorphism \ 
$\phi : G \to {\cal S}(M; \omega) \times \,{\rm ker}\,c^\omega|_G$ \ is defined by \\ 
\hspace{20mm} $\phi(h) = (c^\omega_h, (s(c^\omega_h))^{-1}h)$. \hspace{10mm} 
The inverse is given by \hspace{5mm} $\phi^{-1}(a, f) = s(a)f$. 

(2) Since the topological vector space ${\cal S}(M; \omega)$ strong deformation retracts onto $\{ 0 \}$, 
the conclusion follows from (1). 
\end{proof}


\section*{Appendix} 

This appendix includes the proofs of 
Theorem~\ref{thm_Moser} and Lemma~\ref{lem_collar} stated in Section 2.4.  
These are appropriate parametrized versions of 
Moser's theorem \cite[Theorem 2]{Mos} and 
\cite[Lemma A2]{Mc}, and 
we trace the arguments in \cite{Mos, Mc} together with 
the continuity of related maps with respect to volume forms. 



\begin{proof}[\bf Proof of Theorem~\ref{thm_Moser}] 
For simplicity, let $E_t = E[0, t] (= \partial M \times [0,t])$ and 
$M_t = M - \partial M \times [0,t)$ for $t \in [0,1]$. 
Consider the subspace 
of ${\cal V}^+(M) \times {\cal V}^+(M) \times (0, 1/2)$ 
defined by 
\[ {\cal W}^+(M) = \Big\{ (\sigma, \tau, \delta) \in {\cal V}^+(M) \times {\cal V}^+(M) \times (0, 1/2) \, \Big| \, \sigma(M) = \tau(M),  \sigma = \tau \text{ on } E_{2\delta}\Big\}. \] 
We show that there exists a continuous map $\psi : {\cal W}^+(M) \to {\cal D}_\partial(M)_1$ such that 
$\psi = \psi(\sigma, \tau, \delta)$ satisfies the following conditions: 
\[ \mbox{\rm (1) $\psi_\ast \sigma = \tau$, \ \ 
(2) $\psi = id_M$ on $E_\delta$, \ \ 
(3) if $\sigma = \tau$, then $\psi = id_M$.} \] 
Then the required map $\phi : P \to {\cal D}_\partial(M)_1$ is defined by 
$$\phi(p) = \psi(\mu(p), \nu(p), \e(p)) \ (p \in P).$$

A construction of the map $\psi$ is described in \cite{Mos} and \cite[Ch I \S4]{BT}. 
As for notations, for a smooth orientable $C^\infty$ $n$-manifold $N$ possibly with corners 
and a compact subset $C$ of $N$, let 
$$\Omega^k(N; C) = \big\{ \lambda \in \Omega^k(N) \mid {\rm supp}\,\lambda  \subset C \big\} \ \ (k \geq 0) \hspace{5mm} 
\text{and} \hspace{5mm}   
\Omega^n(N; C)_0 = 
\left\{ \omega \in \Omega^n(N; C) \, \Big| \, \int_N \omega = 0 \right\}.$$ 

The arguments in \cite[Ch I \S4]{BT} (The Poincar\'e Lemma for compactly supported cohomology) yields a continuous map 
\[ \mbox{ 
$\eta_n : \Omega^n({\Bbb R}^n; [-1,1]^n)_0 \to \Omega^{n-1}({\Bbb R}^n; [-1,1]^n)$ \ \ 
such that \ \ $d\hspace{0.2mm}\eta_n(\omega) = \omega$ \ \ 
$(\omega \in  \Omega^n({\Bbb R}^n; [-1,1]^n)_0)$.} \] 
In fact, 
for $\omega \in \Omega^n({\Bbb R}^n; [-1,1]^n)_0$ 
we have 
\[ \omega 
= (-1)^{n-1}dK \omega + (\pi_\ast \omega) \wedge e, 
\] 
where $K\omega \in \Omega^{n-1}({\Bbb R}^n; [-1,1]^n)$ and 
\vspace{1mm} 
$\pi_\ast \omega \in \Omega^{n-1}({\Bbb R}^{n-1}; [-1,1]^{n-1})_0$ are defined as follows. 
We set $(x,t) = ((x_1, \cdots, x_{n-1}), x_n)$. 
Choose any 1-form $e = e(t) dt \in \Omega^1({\Bbb R}; [-1,1])$  such that $\ds \int_{-\infty}^\infty e(s)\,ds = 1$ and put 
$$\ds A(t) = \int_{-\infty}^t e(s)\,ds \in C^\infty({\Bbb R}).$$ 
If $\omega = f \, dx_1 \wedge \cdots \wedge dx_n$ \ \ $(f \in C^\infty({\Bbb R}^n)$, ${\rm supp}\,f \subset [-1, 1]^n)$, then 
\begin{itemize}
\item[] 
$K\omega =  \overline{f} \,dx_1 \wedge \cdots \wedge dx_{n-1}$, \ \ where \ \ 
$\ds \overline{f}(x, t) = \int_{-\infty}^t f(x,s)\,ds - A(t) \int_{-\infty}^\infty f(x,s)\,ds$, \ \ and 
\vskip 1mm 
\item[] $\ds \pi_\ast \omega = \left(\int_{-\infty}^\infty f(x,s)\,ds \right) dx_1 \wedge \cdots \wedge dx_{n-1}$. 
\end{itemize} 

Thus, if the map $\eta_{n-1}$ is already defined, then 
$\pi_\ast \omega = d\hspace{0.2mm}\eta_{n-1}(\pi_\ast \omega)$ and 
the map $\eta_n$ is recursively defined by 
\[ \eta_n(\omega) = (-1)^{n-1}K \omega + \eta_{n-1}(\pi_\ast \omega) \wedge e. \] 

To uniformize the local statement in the preceeding paragraph, the Mayer - Vietoris argument is applied with using a good cover and a partition of unity (cf. \cite[Ch I \S5]{BT}). 
Take any open cover ${\cal U} = \{ U_j \}_{j=1}^m$ of $M_{1/2}$ in ${\rm Int}\, M_{1/3}$
such that ${\cal U}$ is a good cover of $\cup_j U_j$ and $U_j \cap M_{1/2} \neq \emptyset$ $(j=1, \cdots, m)$. 
Then, by the proof of \cite[Lemma 1]{Mos} 
we can find compact subsets $K_j$ in $U_j$ $(j=1, \cdots, m)$ and 
continuous maps  
$$\Omega^n(M; M_{1/2})_0 \ni \omega \longmapsto 
\omega_j \in \Omega^n(M; K_j)_0 \quad (j=1, \cdots, m)$$ 
such that $\omega = \sum_{j=1}^m \omega_j$. 
Since there exists a diffeomorphism from $U_j$ onto ${\Bbb R}^n$ which maps $K_j$ into $[-1, 1]^n$, we can apply the result in the preceeding paragraph to each $U_j$ and sum up to obtain a continuous map 
\[ \eta : \Omega^n(M; M_{1/2})_0 \to \Omega^{n-1}(M; M_{1/3}) \ \ 
\text{such that} \ \ 
d\eta(\omega) = \omega \ \ (\omega \in  \Omega^n(M; M_{1/2})_0). \] 
Replacing $\eta(\omega)$ by $\eta(\omega) - \eta(0)$, 
we may assume that $\eta(0) = 0$. 

Now we can apply the Moser's argument \cite[Theorem 2]{Mos}. 
Define the subspace 
${\cal W}^+(M)_{1/4}$ of \ ${\cal V}^+(M) \times {\cal V}^+(M)$ \ by 
\[ {\cal W}^+(M)_{1/4} = \Big\{ (\sigma, \tau) \in {\cal V}^+(M) \times {\cal V}^+(M) \, \Big| \, \sigma(M) = \tau(M), \ \sigma = \tau \text{ on } E_{1/2}\Big\}. \]
For each $(\sigma, \tau) \in {\cal W}^+(M)_{1/4}$ we have a smooth family of volume forms $\tau_t = \tau + t(\sigma - \tau) \in {\cal V}^+(M)$ $(t \in [0,1])$. 
Since $\sigma - \tau \in \Omega^n(M; M_{1/2})_0$, it follows that 
$\dot{\tau}_t = \sigma - \tau = d\eta(\sigma - \tau)$.
The linear equation 
\begin{equation}\tag{$\ast$}
\iota(X_t) \tau_t = - \eta(\sigma - \tau) 
\end{equation}
uniquely determines a smooth family of vector fields $X_t$ $(t \in [0,1])$ on $M$, 
which induces the associated smooth family of diffeomorphisms 
$\psi_t = \psi_t(\sigma, \tau) \in {\cal D}(M)$ such that $\psi_0 = id_M$. 
The equation $(\ast)$ implies that 
\begin{itemize}
\item[(a)] $\ds \frac{d}{dt} \psi_t^\ast \tau_t 
= \psi_t^\ast\big( \dot{\tau}_t + \iota(X_t) d\tau_t  + d \,\iota(X_t) \tau_t \big) = 0$, \quad hence \quad $\psi_t^\ast \tau_t = \tau_0$ \quad and \quad 
$\psi_1^\ast \sigma = \tau$,
\vskip 2mm 
\item[(b)] if $\sigma = \tau$, then $X_t \equiv 0$ and so $\psi_t \equiv id_M$.
\end{itemize}
Since $\eta(\sigma - \tau) \in \Omega^{n-1}(M; M_{1/3})$, it follows that 
${\rm supp}\,X_t \subset M_{1/3}$ and hence ${\rm supp}\,\psi_t \subset M_{1/3}$. 
The coordinate representation of $X_t$ in \cite[p 293]{Mos} shows that 
the correspondence $(\sigma, \tau) \mapsto (X_t)_{t \in [0,1]}$ is 
continuous under the compact-open $C^\infty$-topology. 
Thus, the correspondence $(\sigma, \tau) \mapsto (\psi_t)_{t \in [0,1]}$ is also continuous under the compact-open $C^\infty$-topology and the map 
$$\psi' : {\cal W}^+(M)_{1/4} \to {\cal D}_{E_{1/3}}(M)_1, \quad \psi' = \psi_1(\sigma, \tau)^{-1}$$ 
satisfies the conditions : \quad 
(c) \ $\psi'_\ast \sigma = \tau$ \quad and \quad  
(d) \ if $\sigma = \tau$, then $\psi' = id_M$.

Finally we incorporate the $(0, 1/2)$-factor to ${\cal W}^+(M)_{1/4}$. 
By a suitable smooth sliding in the collar $E$ we can construct a map 
$$\chi : (0, 1/2) \to {\cal D}_{\partial \cup M_1}(M)_1$$  
such that \quad 
(e) $\chi_\delta = id$ \ on \ $E_{\delta/2}$ \quad and 
\quad (f) $\chi_\delta(E_\delta) =  E_{1/3}$, \ \ $\chi_\delta(E_{2\delta}) =  E_{1/2}$.

Then, the desired map $\psi : {\cal W}^+(M) \to {\cal D}_\partial(M)_1$ is defined by 
$$\psi(\sigma, \tau, \delta) = \chi_\delta^{-1} \psi'((\chi_\delta)_\ast \sigma, (\chi_\delta)_\ast \tau)\chi_\delta.$$ 
This completes the proof. 
\end{proof} 

\begin{proof}[\bf Proof of Lemma~\ref{lem_collar}] 
Since the maps $\mu$, $\nu$ factor through the projections $\pi_1$, $\pi_2$
from the product 
$P_0 = {\cal V}^+(M) \times {\cal V}^+(M)$ onto the 1st and 2nd factor, 
it suffices to consider the case where $(P, \mu, \nu) = (P_0, \pi_1, \pi_2)$. 

(1) We may assume that 
$(M, E, S) = (S \times {\Bbb R}, S \times [-1,1], S \times \{ 0 \})$, 
where $S$ is an oriented $C^\infty$ $(n-1)$-manifold possibly with boundary and 
the orientation of $M$ is induced from that of $S$ and 
the standard orientation of ${\Bbb R}$. 
We fix $\omega_0 \in {\cal V}^+(S)$ and set 
$\omega = \omega_0 \wedge dt \in {\cal V}^+(M)$. 

We shall construct maps 
\[\mbox{
$\psi : {\cal V}^+(S) \to {\cal D}_{S \cup (M - E_U)}(M)_1$ \quad 
and \quad $\delta : {\cal V}^+(S) \to C^\infty(S, (0,1))$
} \]
such that for $\mu \in {\cal V}^+(S)$ 
\begin{itemize} 
\item[(i)$'$\,] $({\psi_\mu})_\ast \mu = \omega$ \ on \ $E_K[-\delta_\mu, \delta_\mu]$, 
\item[(ii)$'$] for each $x \in S$, \ (a) \ $\psi_\mu(E_x^\pm) = E_x^\pm$ \ \ and \ \
\begin{itemize} 
\item[(b)] $(\psi_\mu)|_{E_x^\pm}$ depends only on $\mu|_{E_x^\pm}$ \ 
(i.e., if $\mu, \nu \in {\cal V}^+(M)$ and 
$\mu = \nu$ on $E_x^\pm$, then $\psi_\mu = \psi_\nu$ on $E_x^\pm$). 
\end{itemize} 
\end{itemize} 

Then the required map $\phi : P_0 \to {\cal D}_{S \cup (M - E_U)}(M)_1$ is defined  by $\phi_{(\mu, \nu)} = \psi_\nu^{-1}\psi_\mu$.
Since 
\[ \mbox{$({\psi_\mu})_\ast \mu = \omega = ({\psi_\nu})_\ast \nu$ \ \ on \ \ 
$E_K[-\delta_{(\mu, \nu)}, \delta_{(\mu, \nu)}]$ \quad 
$(\delta_{(\mu, \nu)} = \min\,\{ \delta_\mu, \delta_\nu \})$,} \]
if we choose a map $\varepsilon : P_0 \to C^\infty(S, (0,1))$ 
so that $\varepsilon_{(\mu, \nu)}$ 
is sufficiently small for each $(\mu, \nu) \in P_0$, then 
the maps $\phi$ and $\e$ satisfies the required conditions. 

The construction of the map $\psi$ is described in \cite[Lemma A2]{Mc}. 
The key observation is that for each $f \in C^\infty(M)$ the associated map 
\[ \mbox{$\phi_f : M \to M$, \ $\phi_f(x,t) = (x, f(x,t))$} \] 
has the property that $(\phi_f)^\ast \omega = f_t \,\omega$, where $f_t = \frac{\partial f}{\partial t\,}$. 

Consider the following subspaces of $C^\infty(M)$: 
\[ C^\infty(M)_{>0} = \{ f \in C^\infty(M) \mid f > 0 \}, \quad 
C^\infty(M)' = \{ f \in C^\infty(M) \mid 
f(x,0) = 0, \ f_t > 0 \}. \]


Take three functions $\alpha \in C^\infty({\Bbb R}, (-1/2, 1/2))$, 
$\lambda \in C^\infty({\Bbb R}, [0,1])$ and 
$\rho \in C^\infty(S, [0,1])$ such that 
\begin{itemize}
\item[(a)] $\alpha(t) = t$ $(|t| \leq 1/3)$, \ $\alpha_t > 0$, 
\item[(b)] $\lambda(t) = \lambda(-t)$, \ $\lambda(t) = 0$ $(|t| \leq 1/2)$, \ 
$\lambda(t) = 1$ $(|t| \geq 1)$, \ $\lambda_t(t) \geq 0$ $(t \geq 0)$ \ and \
\item[(c)] $\rho(x) = 0$ $(x \in V)$, \ $\rho(x) = 1$ on $(x \in S - U)$, 
where $V$ is an open subset of $S$ with $K \subset V \subset U$. 
\end{itemize} 
Choose a continuous map $\gamma : C^\infty(M)' \to C^\infty(S, (0,1/2))$
such that 
\[ f(E[-\gamma_f, \gamma_f]) \subset [-1/3, 1/3] \quad (f \in C^\infty(M)'). \]

For each $f \in C^\infty(M)'$ we define 
$\widetilde{f} \in C^\infty(M)'$ by 
\[ \widetilde{f}(x,t) = (1 - \rho(x)) g(x,t) + \rho(x) \cdot t, \qquad g(x,t) = (1 - \lambda(t))\alpha(f(x,t)) + \lambda(t)\cdot t \qquad 
((x,t) \in M). \]  
Since $\widetilde{f}(x,t) = t$ $((x,t) \in S \cup (M - E_U))$ and 
$\widetilde{f} = f$ on $E_V[-\gamma_f, \gamma_f]$, it follows that 
\[ \mbox{
(a) \ $\phi_{\widetilde{f}} \in {\cal D}_{S \cup (M - E_U)}(M)_1$ \quad and \quad 
(b) \ $(\phi_{\widetilde{f}})^\ast \omega = \widetilde{f}_t \,\omega = f_t \,\omega$ \ on \ $E_V[-\gamma_f, \gamma_f]$.} \]  

Each $\mu \in {\cal V}^+(M)$ uniquely 
determines $h_\mu \in C^\infty(M)_{>0}$ such that  
$\mu = h_\mu \omega$. 
Define $f_\mu \in C^\infty(M)'$ by 
\[ f_\mu(x,t) = \int_0^t h_\mu(x,s)\,ds. \] 
Since $(f_\mu)_t = h_\mu$, it is seen that 
\[ \mbox{
$(\phi_{\widetilde{f_\mu}})^\ast \omega = h_\mu \,\omega = \mu$ \ on \ $E_V[-\gamma_{f_\mu}, \gamma_{f_\mu}]$.} \] 

Finally we define the map 
$\psi : {\cal V}^+(S) \to {\cal D}_{S \cup (M - E_U)}(M)_1$ by $\psi_\mu = \phi_{\widetilde{f_\mu}}$.  
The condition (ii)$'$ follows from the construction of the map $\psi$. 
If we choose a map $\delta : {\cal V}^+(S) \to C^\infty(S, (0,1))$ 
so that $\delta_\mu$ 
is sufficiently small for each $\mu \in {\cal V}^+(S)$, then 
the maps $\psi$ and $\delta$ satisfy the condition (i)$'$. 

(2) We may assume that 
$(M, E, \partial M) = (S \times [0, \infty), S \times [0,1], S \times \{ 0 \})$, 
where $S$ is an oriented $C^\infty$ $(n-1)$-manifold possibly with boundary. 
The assertion (2) follows from the same argument as in (1) by ignoring the $(-\infty, 0]$-side. 
\end{proof}


\end{document}